\documentstyle[12pt]{article}
\title{{\bf Algorithms for $D$-modules ---restriction, tensor product, 
localization, and local cohomology groups}}
\author{Toshinori Oaku and Nobuki Takayama}
\date{May 2, 1998}

\newcommand{\ord}{\mbox{{\rm ord}}}
\newcommand{\lexp}{{\rm lexp}}
\newcommand{\lcoef}{{\rm lcoef}}

\newcommand{\mod}{\quad \mbox{{\rm mod}}\quad }
\newcommand{\gr}{\mbox{{\rm gr}}}

\newcommand{\Ext}{\mbox{{\rm Ext}}}

\newcommand{\Bsc}{{\cal B}}
\newcommand{\Dsc}{{\cal D}}

\newcommand{\Hsc}{{\cal H}}
\newcommand{\Isc}{{\cal I}}
\newcommand{\Jsc}{{\cal J}}
\newcommand{\Ksc}{{\cal K}}
\newcommand{\Lsc}{{\cal L}}
\newcommand{\Msc}{{\cal M}}
\newcommand{\Nsc}{{\cal N}}
\newcommand{\Osc}{{\cal O}}

\newcommand{\Homsc}{{\cal H}om\,}
\newcommand{\Extsc}{{\cal E}xt\,}
\newcommand{\Torsc}{{\cal T}or\,}
\newcommand{\G}{{G}}
\newcommand{\N}{{\bf N}}
\newcommand{\Z}{{\bf Z}}
\newcommand{\Q}{{\bf Q}}

\newcommand{\C}{{\bf C}}
\newcommand{\Set}[1]{\{1,\dots,#1\}}
\newcommand{\mvec}{{\bf m}}
\newcommand{\nvec}{{\bf n}}
\newcommand{\init}{\mbox{{\rm in}}}
\newcommand{\lp}{\mbox{{\rm lp}}}
\newcommand{\an}[2]{{\cal #1}^{\rm an}_{#2}}
\def\pd#1{ \partial_{#1} }

\def\upmap#1{\rule[-12pt]{0pt}{30pt}\Big\uparrow\rlap
  {$\vcenter{\hbox{$\scriptstyle#1$}}$}}
\newenvironment{proof}{\par\noindent Proof:\ }{${\tt [}\kern-0.2mm{\tt ]}$}
\newtheorem{theorem}{Theorem}[section]
\newtheorem{proposition}[theorem]{Proposition}
\newtheorem{lemma}[theorem]{Lemma}

\newtheorem{definition}[theorem]{Definition}
\newtheorem{remark}[theorem]{Remark}
\newtheorem{example}[theorem]{Example}
\newtheorem{algorithm}[theorem]{Algorithm}
\begin{document}

\maketitle

\section{Introduction}

The purpose of this paper is to describe algorithms for 
computing various functors for algebraic $\Dsc$-modules,   
i.e.\  systems of linear partial differential equations 
with polynomial coefficients.
The algorithms enable us to perform actual computations (with
limitation caused by the complexity) by using e.g.\ a program  
{\tt kan} \cite{Kan-sm1} developed by the second author,  
as well as to establish theoretical 
computability of some fundamental functors in the $\Dsc$-module theory.   

Let $K$ be an algebraically closed field of characteristic zero and let 
$X$ be the affine space $K^n$ with a positive integer $n$.  
We denote by $\Osc_X$ and $\Dsc_X$ the sheaves on $X$ of 
rings of regular functions and of algebraic linear differential operators 
respectively (cf.\ Bernstein \cite{BernsteinNotes}, 
Bj\"ork \cite{BjorkBook}, Borel et al.\ \cite{Borel}).
Let $\Msc$ and $\Nsc$ be coherent left $\Dsc_X$-modules.

Various functors are defined for (especially for holonomic) $\Dsc$-modules
and play the fundamental 
role (see \cite{BernsteinNotes},\cite{BjorkBook},\cite{Borel} and also e.g., 
Kashiwara \cite{HolonomicII}, \cite{KashiwaraBook}, 
Mebkhout \cite{MebkhoutBook} for their analytic counterparts). 
Among such functors, we are concerned with the following:
\begin{enumerate}
\item
The cohomology groups of the restriction 
$\Msc_Y^\bullet := \Osc_Y\otimes^{L}_{\Osc_X}\Msc$ of $\Msc$ to $Y$ 
as left $\Dsc_Y$-modules, 
where $Y$ is a non-singular subvariety of $X$ and 
$\otimes^{L}$ denotes the left derived functor 
(cf.\ \cite{HartshorneNotes}) of the tensor product.
\item
The cohomology groups $\Ext^i_{\Dsc_X}(\Msc, K[[x_1,\dots,x_n]])$ 
with coefficients in the formal power series solutions of $\Msc$,  
which equal to those with coefficients 
in the convergent power series solutions if $\Msc$ 
is regular holonomic (cf.\ Kashiwara-Kawai \cite{KKholonomic3}).
\item
The tensor product $\Msc\otimes_{\Osc_X}\Nsc$ and, more generally, the torsion
groups $\Torsc^i_{\Osc_X}(\Msc,\Nsc)$, as left $\Dsc_X$-modules.
\item
The localization 
$\Msc[f^{-1}] := \Osc_X[f^{-1}]\otimes_{\Osc_X}\Msc$ of $\Msc$ as a left
$\Dsc_X$-module, where
$f \in K[x_1,\dots,x_n]$ is an arbitrary non-constant polynomial.  
\item
The (algebraic) local cohomology groups $\Hsc^i_{[Y]}(\Msc)$ with 
support $Y$ as left $\Dsc_X$-modules, 
where $Y$ is an arbitrary algebraic set of $X$. 
\end{enumerate} 
It was proved by Kashiwara \cite{HolonomicII} that these are all holonomic
systems (the second one is a finite dimensional vector space)
if so are $\Msc$ and $\Nsc$.  

Let us remark that if $K= \C$ and $\Msc$ is Fuchsian along $Y$ in the sense 
of Laurent and Moteiro-Fernandes \cite{LM}, 
which is the case if $\Msc$ is regular holonomic 
in the sense of \cite{KKholonomic3}, then there exists an isomorphism 
$$
R\Homsc_{\Dsc_X}(\Msc,\Osc^{\rm an}_X)|_Y \simeq
R\Homsc_{\Dsc_Y}(\Msc^\bullet_Y,\Osc^{\rm an}_Y)
$$
in the derived category of sheaves of $\C$-vector spaces; 
here $\Osc^{\rm an}_X$ and $\Osc^{\rm an}_Y$ denotes the sheaves of 
holomorphic functions on $X$ and on $Y$ respectively, 
and $R\Homsc$ the right derived functor of $\Homsc$.
Thus roughly speaking, $\Msc^\bullet_Y$ corresponds to the system of 
partial differential equations which the solutions of $\Msc$ 
restricted to $Y$ satisfy.   
Similarly, $\Msc\otimes_{\Osc_X}\Nsc$ corresponds to the system which 
the product of solutions of $\Msc$ and of $\Nsc$ satisfies.   

As was observed by Galligo \cite{Galligo} and was developed by several authors 
(e.g.\ \cite{Castro}, \cite{TakayamaJJAM}, \cite{TakayamaISSAC90}, 
\cite{OakuChar}, \cite{OakuFuchs}, \cite{OakuBfunction}, \cite{ACG}, 
\cite{Sturmfels-Takayama})
the notion of Gr\"obner basis and the Buchberger algorithm \cite{Buchberger}
are essential in the algorithmic study of $\Dsc$-modules as well as 
in computational algebraic geometry 
(cf.\ \cite{CoxBook}, \cite{EisenbudBook}). 
By using Gr\"obner bases for the Weyl algebra, 
we give algorithms for computing the objects listed above 
under some conditions on $\Msc$ and $\Nsc$, 
which are certainly satisfied if $\Msc$ and $\Nsc$ are holonomic.
These algorithms also apply to the analytic counterparts of these functors 
as long as the input $\Dsc$-module is defined algebraically. 

We first give an algorithm for the restriction 
(Algorithm \ref{algorithm:restriction})
when $Y$ is a linear subvariety of arbitrary codimension 
under the condition that $\Msc$ is specializable along $Y$, 
which is the case with an arbitrary holonomic $\Dsc_X$-module $\Msc$.  
Here $\Msc$ is specializable along $Y$ by definition 
if and only if there exists a nonzero $b$-function, 
or the indicial polynomial of $\Msc$ along $Y$. 
We also give an algorithm to compute the $b$-function 
(Algorithm \ref{algorithm:b-function}).  

Our method consists in computing a free resolution of $\Msc$ that is 
adapted to the so-called $V$-filtration associated with $Y$. 
Such a free resolution tensored with 
$\Dsc_{Y\rightarrow X} := (\Dsc_X)_Y^\bullet$ gives $\Msc_Y^\bullet$, 
but it is not a complex of coherent $\Dsc_Y$-modules in general.   
Then we use information on the integral roots of the $b$-function 
to truncate the complex and obtain a complex of finitely generated 
free $\Dsc_Y$-modules. 
The first author gave in \cite{OakuAdvance} an algorithm
for the case where $Y$ is of codimension one 
without using free resolution. 

This algorithm for the restriction also solves the other problems 
by virtue of some isomorphisms provided by the $\Dsc$-module theory, 
especially those described in \cite{HolonomicII}. 
See Algorithm \ref{algorithm:tensor} for the tensor product, 
Algorithm \ref{algorithm:localization} for the localization, 
and Algorithm \ref{algorithm:localCohomology} for the algebraic 
local cohomology groups.
Finally the computation of the restriction for the general case where 
$Y$ is not necessarily linear reduces to that of local cohomology  
through the so-called Kashiwara equivalence \cite{HolonomicII}, 
which claims the equivalence of the category of coherent $\Dsc_Y$-modules 
and that of coherent $\Dsc_X$-modules supported by $Y$. 

Algorithms for the local cohomology groups 
have been given in \cite{OakuAdvance} 
when $Y$ is of codimension one, and by Walther \cite{Walther} under the 
assumption that $\Msc$ is saturated with respect to $Y$.
An algorithm for the localization has been given in \cite{OakuAdvance} 
under the condition that $\Msc$ is $f$-saturated.  
These assumptions are removed in the present paper. 
As another application of the restriction algorithm, 
we can also obtain an algorithm for integration of a module over the 
Weyl algebra, which enables us to compute the de Rham cohomology groups
of some algebraic varieties. See \cite{OTdeRham} for details.
 
Finally we discuss how to get
the free resolution mentioned above.  
For that purpose, we apply Schreyer's 
method for free resolution in the polynomial ring 
(see e.g.\ \cite{EisenbudBook})  
to the ring of differential operators.  
In doing so, we need some modification because of 
the non-commutativity and 
the fact that the term order we use is not a well-order.
We have two methods to cope with this difficulty:
one is the homogenization with respect to the $V$-filtration 
by the first author (\cite{OakuFuchs},\cite{OakuBfunction}); 
the other is what we call the homogenized Weyl algebra  
which was introduced and implemented by the second author 
in the 2nd version of {\tt kan/sm1} \cite{Kan-sm1} that was released in 1994, 
but has not been published in the literature.    
A similar method was employed by Assi et al.\ \cite{ACG} independently 
and applied to the computation of the slopes of a $\Dsc$-module.  

We have implemented the algorithms by using {\tt kan/sm1} \cite{Kan-sm1} 
for computations of Gr\"obner bases and free resolutions in the Weyl algebra,
and Risa/Asir \cite{Risa} for factorization and primary decomposition 
in the polynomial ring. 

\section{V-filtration and free resolution}
\label{section:V-filtration}
\setcounter{equation}{0}

Let $K$ be an algebraically closed field of characteristic zero.  
We fix  positive integers $d$ and $n$. 
Let $X$ be the affine space $K^{d+n}$ with the 
coordinate system $(t,x) = (t_1,...,t_d,x_1,\dots,x_n)$.  
We denote by $\pd{t} = (\pd{t_1},...,\pd{t_d})$ 
and $\pd{x} = (\pd{x_1},\dots,\pd{x_n})$ 
the corresponding derivations with 
$\pd{x_i} = \partial/\partial x_i$, $\pd{t_j}=\partial/\partial t_j$.  
We use the notation 
$x^\alpha := x_1^{\alpha_1}\cdots x_n^{\alpha_n}$, 
$\pd{x}^\beta := \pd{x_1}^{\beta_1}\cdots \pd{x_n}^{\beta_n}$,  
$t^\mu := t_1^{\mu_1}\cdots t_d^{\mu_d}$, 
$\pd{t}^\nu := \pd{t_1}^{\nu_1}\cdots\pd{t_n}^{\nu_n}$ 
for $\alpha = (\alpha_1,\dots,\alpha_n)$, 
$\beta = (\beta_1,\dots,\beta_n) \in \N^n$ and
$\mu = (\mu_1,\dots,\mu_d), \nu = (\nu_1,...,\nu_d) \in \N^d$,
where we put $\N := \{0,1,2,\dots \}$.   
We also use the notation $|\alpha| := \alpha_1 + \dots + \alpha_n$.  

Let $Y$ be the $d$-codimensional linear subvariety of $X$ given by 
$Y := \{(t,x) \in X \mid t=0\}$.  
Let $\Osc_X$ and $\Osc_Y$ be the sheaves of 
regular functions on $X$ and on $Y$ respectively.  
We denote by $\Dsc_X$ and $\Dsc_Y$ the sheaves of 
rings of algebraic linear differential operators on $X$ and on $Y$ 
respectively.  

Let $\Msc$ be a coherent left $\Dsc_X$-module on $X$.
Then the set of the global sections $M := \Gamma(X,\Msc)$ is 
a finitely generated left module over the Weyl algebra 
$A_{d+n} := \Gamma(X,\Dsc_X)$.
Conversely, for a finitely generated left $A_{d+n}$-module $M$,
its sheafification $\Msc := \Dsc_X\otimes_{A_{d+n}}M$ is a coherent 
$\Dsc_X$-module.  More precisely, this correspondence gives an equivalence
between the category of finitely generated $A_{d+n}$-modules and that of
coherent $\Dsc_X$-modules (cf.\ \cite{BernsteinNotes},\cite{Borel}).
Hence we could work only in the first category.  
However, as to e.g., the restriction functor, it would be preferable to work 
in the latter category since a coherent 
$\Dsc_X$-module can be specializable along some Zariski open subset of 
$Y$ but not along whole $Y$ (cf.\ Section \ref{section:b-function}).  
In any case, actual computations are done for modules over the Weyl
algebra.  

In the sequel, we define the notion of free resolution adapted to the
$V$-filtration. 
Let $\Dsc_X|_Y$ be the sheaf theoretic restriction
of $\Dsc_X$ to $Y$.  
Let $\Jsc_Y := \Osc_Xt_1 + \dots + \Osc_Xt_d$ be the defining ideal of $Y$.  
Then for each integer $k$ we put
\begin{eqnarray*}
F_Y^k(\Dsc_X) &:=& \{ P \in \Dsc_X|_Y \mid P (\Jsc_Y)^j 
\in (\Jsc_Y)^{j-k} \quad\mbox{for any $j\geq k$} \} \\
&=& 
\{ P = \sum_{|\nu|\leq l}\sum_{|\beta|\leq m}
 a_{\nu\beta}(t,x)\pd{t}^\nu\pd{x}^\beta  \mid 
l,m \in \N,\,\, a_{\nu\beta}(t,x) \in \Jsc_Y^{|\nu|-k}\}
\end{eqnarray*}
with the convention $\Jsc_Y^j = \Osc_X$ for $j \leq 0$.
This is called the V-filtration attached to $Y$ 
(cf.\ \cite{KashiwaraVanishing},\cite{LS}).  
More generally, given an $r$-vector 
$\mvec := (m_1,...,m_r)$ of integers, we put  
$$
F_Y^k[\mvec](\Dsc_X^r) := \bigoplus_{i=1}^r F_Y^{k-m_i}(\Dsc_X)e_i,
$$
where $e_1,\dots,e_r$ are the canonical generators of $\Dsc_X^r$.    
We may assume that $\Msc$ has a presentation
$\Msc = \Dsc_X^r/\Nsc$ on $X$, 
where $\Nsc$ is a coherent left $\Dsc_X$-submodule of 
$\Dsc_X^r$.
In fact, $M := \Gamma(X,\Msc)$ can be written in the form 
$M = A_{d+n}^r/N$ with an integer $r$ and an $A_{d+n}$-submodule $N$ 
of $A_{d+n}^r$.  Then $\Nsc := \Dsc_X\otimes_{A_{d+n}}N$ satisfies 
the above property. 
Let $u_i$ be the residue class of $e_i \in \Dsc_X^r$ in $\Msc$.  
Then for $\mvec \in \Z^r$, we put
\begin{eqnarray*}
F_Y^k[\mvec](\Nsc) &:=& \Nsc \cap F_Y^k[\mvec](\Dsc_X^r), \\
F_Y^k[\mvec](\Msc) &:=&  
 F_Y^{k-m_1}(\Dsc_X)u_1 + \cdots + F_Y^{k-m_r}(\Dsc_X)u_r
\end{eqnarray*}
for each integer $k \in \Z$. 
The graded ring and modules associated 
with these filtrations are defined by
\begin{eqnarray*}
\gr_Y(\Dsc_X) &:=& 
  \bigoplus_{k\in\Z}F_Y^k(\Dsc_X)/F_Y^{k-1}(\Dsc_X), \\
\gr_Y[\mvec](\Dsc_X^r) &:=& 
  \bigoplus_{k\in\Z}F_Y^k[\mvec](\Dsc_X^r)/F_Y^{k-1}[\mvec](\Dsc_X^r), \\
\gr_Y[\mvec](\Nsc) &:=& \bigoplus_{k\in\Z} 
  F_Y^k[\mvec](\Nsc)/F_Y^{k-1}[\mvec](\Nsc). \\
\end{eqnarray*}
Then $\gr_Y[\mvec](\Nsc)$ and $\gr_Y[\mvec](\Msc)$ are 
coherent left $\gr(\Dsc_X)$-modules. 
If $\mvec$ is the zero vector, we shall omit the notation $[\mvec]$.  

For a nonzero section $P$ of $\Dsc_X^r|_Y$, let $k = \ord_Y[\mvec](P)$ be 
the minimum $k \in \Z$ such that $P \in F_Y^k[\mvec](\Dsc_X^r)$.   
(We put $\ord_Y[\mvec](0) := -\infty$.)
Then let $\sigma_Y[\mvec](P)$ be the residue class of $P$ in 
$$
\gr_Y^k[\mvec](\Dsc_X^r) 
:= F_Y^{k}[\mvec](\Dsc_X^r)/F_Y^{k-1}[\mvec](\Dsc_X^r).  
$$

\begin{definition} \label{def:resolution}
\rm
Let $\Msc = \Dsc_X^r/\Nsc$ be as above.   
Let us consider a free resolution 
\begin{equation} \label{eq:resolution}
\Dsc_X^{r_l} \stackrel{\psi_l}{\longrightarrow} 
\Dsc_X^{r_{l-1}} 
\stackrel{\psi_{l-1}}{\longrightarrow} 
\cdots 
\stackrel{\psi_2}{\longrightarrow} \Dsc_X^{r_1}
\stackrel{\psi_1}{\longrightarrow} \Dsc_X^{r_0}
\stackrel{\varphi}{\longrightarrow}\Msc \longrightarrow 0
\end{equation}
of $\Msc$, 
where $\psi_i$ are homomorphisms of left $\Dsc_X$-modules, 
and $\varphi$ is defined by $\varphi(e_i) = u_i$ for $i=1,...,r_0$ 
with $r_0 = r$.  
This free resolution is said to be {\em adapted} (or {\em strict} 
with respect) to the V-filtration $F_Y[\mvec]$
if and only if there exist vectors $\mvec_1 \in \Z^{r_1}$,$...$,
$\mvec_l \in \Z^{r_l}$ such that 
$$
\psi_{j+1}(F_Y^k[\mvec_{j+1}](\Dsc_X^{r_{j+1}}))
\subset F_Y^k[\mvec_j](\Dsc_X^{r_j})
$$ 
holds for $j=0,1,\dots,l-1$ with $\mvec_0 = \mvec$ and that 
$$
F_Y^k[\mvec_l](\Dsc_X^{r_l}) \stackrel{\psi_l}{\longrightarrow} 
\cdots 
\stackrel{\psi_2}{\longrightarrow} F_Y^k[\mvec_1](\Dsc_X^{r_1})
\stackrel{\psi_1}{\longrightarrow} F_Y^k[\mvec_0](\Dsc_X^{r_0})
\stackrel{\varphi}{\longrightarrow}
F_Y^k[\mvec](\Msc) \rightarrow 0
$$
is an exact sequence for any $k \in \Z$.
We call $\mvec_1$,...,$\mvec_l$ the {\em shift vectors} 
associated with the free resolution (\ref{eq:resolution}).
\end{definition}

The definition above is general in the sense that it is local and also 
applies to the analytic case (cf.\ Section \ref{section:analytic}).  
However, from the computational view point, working in the Weyl algebra 
would be more convenient: 
Let $A_n$ and $A_{d+n}$ be the Weyl algebras on the $n$ 
variables $x$ and on the $d+n$ variables $(t,x)$ respectively 
with coefficients in $K$ (cf.\ \cite{BjorkBook}).  
Put $L := \N^{2(d+n)} = \N^d \times \N^d \times \N^n \times \N^n$.
An element $P$ of $A_{d+n}^r$ is written in a finite sum
\begin{equation} \label{eq:PinWeyl}
P = \sum_{i=1}^r \sum_{(\mu,\nu,\alpha,\beta) \in L}
       a_{\mu\nu\alpha\beta i} t^\mu\pd{t}^\nu x^\alpha\pd{x}^\beta e_i
\end{equation}
with $a_{\mu\nu\alpha\beta i} \in K$,
$e_1 := (1,0,\dots,0),\dots,e_r := (0,\dots,0,1)$.
Put
\begin{eqnarray*}
F_Y^k(A_{d+n}) &:=& \{ P = \sum_{|\nu|-|\mu|\leq k}\sum_\beta
 a_{\mu\nu\beta}(x)t^\mu\pd{t}^\nu\pd{x}^\beta  \in A_{d+n} \mid 
a_{\mu\nu\beta}(x) \in K[x] \},\\
F_Y^k[\mvec](A_{d+n}^r) &:=& \bigoplus_{i=1}^r F_Y^{k-m_i}(A_{d+n})e_i,\\
F_Y^k[\mvec](M) &:=& F_Y^{k-m_1}(A_{d+n})u_1 + \cdots + F_Y^{k-m_r}(A_{d+n})u_r,
\end{eqnarray*}
where $u_i$ is the residue class of $e_i$ in $M$.  

The following lemma follows immediately from the definition:

\begin{lemma} \label{lemma:AtoD}
Let $p$ be a point of $Y$ and $P$ a germ of $\Dsc_X$ at $p$.  
Then $P$ belongs to $F_Y^k(\Dsc_X)_p$ if and only if there exists 
$a(t,x) \in K[t,x]$ such that $a(p) \neq 0$ and $a(t,x)P \in F_Y^k(A_{d+n})$.
\end{lemma}


By using this lemma and the flatness of $\Dsc_X$ over $A_{d+n}$, we 
can easily get the following:

\begin{proposition} \label{prop:AtoD}
Let 
\begin{equation} \label{eq:A-resolution}
A_{d+n}^{r_l} \stackrel{\psi_l}{\longrightarrow} 
A_{d+n}^{r_{l-1}} 
\stackrel{\psi_{l-1}}{\longrightarrow} 
\cdots 
\stackrel{\psi_2}{\longrightarrow} A_{d+n}^{r_1}
\stackrel{\psi_1}{\longrightarrow} A_{d+n}^{r_0}
\stackrel{\varphi}{\longrightarrow} M \longrightarrow 0
\end{equation}
be a free resolution of $M$ adapted to the $F_Y[\mvec]$-filtration, i.e.,  
$\psi_i$ are homomorphisms of left $A_{d+n}$-modules, 
$\varphi$ is defined by $\varphi(e_i) = u_i$ for $i=1,...,r_0$ 
with $r_0 = r$, and there exist $\mvec_1 \in \Z^{r_1}$,$...$,
$\mvec_l \in \Z^{r_l}$ such that 
$$
\psi_{j+1}(F_Y^k[\mvec_{j+1}](A_{d+n}^{r_{j+1}}))
\subset F_Y^k[\mvec_j](A_{d+n}^{r_j})
$$ 
holds for $j=0,1,\dots,l-1$ with $\mvec_0 = \mvec$ and that 
$$
F_Y^k[\mvec_l](A_{d+n}^{r_l}) \stackrel{\psi_l}{\longrightarrow} 
\cdots 
\stackrel{\psi_1}{\longrightarrow} F_Y^k[\mvec_0](A_{d+n}^{r_0})
\stackrel{\varphi}{\longrightarrow}
F_Y^k[\mvec](M) \rightarrow 0
$$
is an exact sequence for any $k \in \Z$.
Under this assumption, the exact sequence (\ref{eq:A-resolution}) 
tensored by $\Dsc_X$ from the left gives 
a free resolution of $\Msc$ adapted to the $F_Y[\mvec]$-filtration.
\end{proposition}

\section{Gr\"obner bases and free resolution} \label{section:Gbase}
\setcounter{equation}{0}

The purpose of this section is to show that Gr\"obner bases homogenized 
with respect to the V-filtration provide a free resolution 
adapted to the V-filtration.  
An alternative and more efficient method will be described in Section
\ref{section:Schreyer}.  

We fix a natural number $r$ and a vector $\mvec \in \Z^r$.
Let $\prec$ be a well-order (i.e. a linear order) 
on $L \times \Set{r}$ which satisfies
\begin{eqnarray} \label{eq:monomial-order}
& &
\mbox{
$(\tilde\alpha,i) \prec (\tilde\beta,j)$ implies 
$(\tilde\alpha+\tilde\gamma,i) \prec (\tilde\beta+\tilde\gamma,j)$
}\nonumber\\
& & 
\mbox{ 
for any $\tilde\alpha,\tilde\beta,\tilde\gamma \in L$ and $i,j \in \Set{r}$.
}
\end{eqnarray}
Then we define a total order $\prec_F$ on $L \times \{1,...,r\}$ by 
\begin{eqnarray} \label{eq:F-order}
& &
\mbox{
$(\mu,\nu,\alpha,\beta,i) \prec_F (\mu',\nu',\alpha',\beta',j)$ 
if and only if 
}\nonumber \\
& &
\mbox{
$|\nu-\mu|+m_i < |\nu'-\mu'|+m_j$ or else
}\nonumber \\
& &
\mbox{
$|\nu-\mu|+m_i = |\nu'-\mu'|+m_j$,  
$(\mu,\nu,\alpha,\beta,i) \prec (\mu',\nu',\alpha',\beta',j)$.
}
\end{eqnarray}
Let $P$ be a nonzero element of $A_{d+n}^r$ which is written in the form
(\ref{eq:PinWeyl}).  
Then the {\em leading exponent} $\lexp_F(P) \in L \times \Set{r}$
of $P$ with respect to $\prec_F$ is defined as the maximum element of
$\{(\mu,\nu,\alpha,\beta,i) \mid a_{\mu\nu\alpha\beta i} \neq 0\}$ 
in the order $\prec_F$. 
Moreover, for $(\mu,\nu,\alpha,\beta,i) = \lexp_F(P)$, the leading 
coefficient of $P$ is defined by $\lcoef_F(P) := a_{\mu\nu\alpha\beta i}$.  
The set of leading exponents $E_F(N)$ of a subset $N$ of $A_{n+1}^r$ 
is defined by
$$ E_F(N) := \{ \lexp_F(P) \mid P \in N \setminus \{0\} \}.$$ 

\begin{definition} \label{def:Gbase}
\rm
A finite set $\G$ of generators of a left $A_{d+n}$-submodule $N$ 
of $A_{d+n}^r$ is called a Gr\"obner basis of $N$ 
with respect to $\prec_F$ 
(or an $F[\mvec]$-Gr\"obner basis),
if we have 
$$ E_F(N) = \bigcup_{P \in \G} (\lexp(P) + L), $$
where we write 
$$ (\tilde\alpha,i) + L = \{ (\tilde\alpha + \tilde\beta,i) \mid  
\tilde\beta \in L \}$$
for $\tilde\alpha \in L$ and $i \in \Set{r}$.  
\end{definition}

We define an order $\prec_H$ on $\N \times L \times \{1,...,r\}$ by
\begin{eqnarray} \label{eq:H-order}
& & 
\mbox{
$(\lambda,\tilde\alpha,i) \prec_H (\lambda',\tilde\alpha',i)$ if and only if
$\lambda < \lambda'$ or else 
} \nonumber \\
& &
\mbox{
$\lambda = \lambda'$,
$(\lambda,\tilde\alpha,i) \prec (\lambda',\tilde\alpha',j)$,
}
\end{eqnarray}
where $\lambda,\lambda' \in \N, \tilde\alpha,\tilde\alpha'\in L$ and 
$i,j \in \{1,...,r\}$.  
It is easy to see that $\prec_H$ is a well-order and satisfies 

\begin{lemma}
If $|\nu-\mu|+m_i-\lambda = |\nu'-\mu'|+m_j-\lambda'$, 
then we have 
$$
(\lambda,\mu,\nu,\alpha,\beta,i) \prec_H (\lambda',\mu',\nu',\alpha',\beta',j)
$$ 
if and only if 
$$
(\mu,\nu,\alpha,\beta,i) \prec_F (\mu',\nu',\alpha',\beta',j).
$$
\end{lemma}

We introduce an indeterminate $t_0$ which commutes with any element of 
$A_{d+n}$ in order to define the homogenization.

\begin{definition} \label{def:F-homogeneous}
\rm
An element $P$ of $A_{d+n}[t_0]^r$ of the form
$$ P = \sum_{i=1}^r \sum_{\lambda,\mu,\nu,\alpha,\beta} 
a_{\lambda\mu\nu\alpha\beta i}t_0^\lambda 
t^\mu x^\alpha\pd{t}^\nu\pd{x}^\beta e_i $$
is said to be $F[\mvec]$-{\em homogeneous} of order $k$ if 
$a_{\lambda\mu\nu\alpha\beta i} = 0$ whenever $|\nu-\mu|-\lambda+m_i \neq k$.  
\end{definition}

\begin{definition} \label{def:F-homogenization}
\rm
For an element $P$ of $A_{d+n}^r$ of the form (\ref{eq:PinWeyl}), 
put 
$$
k := \min\{|\nu-\mu|+m_i \mid a_{\mu\nu\alpha\beta i} \neq 0 
\mbox{ for some $\alpha,\beta$} \}.
$$  
Then the $F[\mvec]$-{\em homogenization} $h(P)\in A_{d+n}[t_0]^r$ 
of $P$ is defined by
$$ 
h(P) := \sum_{i=1}^r\sum_{\mu,\nu,\alpha,\beta} a_{\mu\nu\alpha\beta i}
t_0^{|\nu-\mu|+m_i-k}t^\mu x^\alpha\pd{t}^\nu\pd{x}^\beta e_i.
$$
Then $h(P)$ is $F[\mvec]$-homogeneous of order $k$.
\end{definition}

When $\mvec$ is the zero vector, we simply say $F$-homogeneous 
instead of $F[\mvec]$-homogeneous.  

\begin{lemma}
If $P \in A_{d+n}[t_0]$ is $F$-homogeneous and 
$Q \in A_{d+n}[t_0]^r$ is $F[\mvec]$-homogeneous, 
then $PQ$ is $F[\mvec]$-homogeneous.
\end{lemma}

\begin{lemma}
For $P_1,\dots,P_k \in A_{d+n}^r$, put $P = P_1 + \cdots + P_k$.  
Then there exist $l, l_1,\dots, l_k \in \N$ so that
$$ t_0^l h(P) = t_0^{l_1} h(P_1) + \cdots + t_0^{l_k}h(P_k).$$
\end{lemma}

Let us define 
$\varpi : \N \times L \times \Set{r} \longrightarrow L \times \Set{r}$ by 
$\varpi(\lambda,\mu,\nu,\alpha,\beta,i) = (\mu,\nu,\alpha,\beta,i)$. 
For a nonzero element $P = P(t_0)$ of $A_{d+n}[t_0]^r$, 
let us denote by 
$\lexp_H(P) \in \N \times L \times \Set{r}$ and $\lcoef_H(P) \in K$ 
the leading exponent and the leading coefficient of $P$ 
with respect to $\prec_H$.  

\begin{lemma}
\begin{enumerate}
\item If $P(t_0) \in A_{d+n}[t_0]^r$ is $F[\mvec]$-homogeneous, then 
we have $\lexp_F(P(1)) = \varpi(\lexp_H(P(t_0)))$.
\item For any $P \in A_{d+n}^r$, 
we have $\lexp_F(P) = \varpi(\lexp_H(h(P)))$.   
\end{enumerate}
\end{lemma}

Since the Buchberger algorithm preserves the $F[\mvec]$-homogeneity, 
we have

\begin{proposition}
Let $N$ be a left $A_{d+n}[t_0]$-submodule of $A_{d+n}[t_0]^r$ 
generated by $F[\mvec]$-homogeneous operators.  
Then there exists a Gr\"obner basis 
with respect to $\prec_H$ of $N$ 
consisting of $F[\mvec]$-homogeneous operators.  
Moreover, such a Gr\"obner basis can be computed by the Buchberger 
algorithm.
\end{proposition}

The following proposition can be easily proved in the same way as 
\cite[Theorem 3.12]{OakuBfunction}

\begin{proposition} \label{prop:F-Gbase}
Let $N$ be a left $A_{d+n}$-submodule of $A_{d+n}^r$ generated by 
$P_1,\dots,P_l \in A_{d+n}^r$.  
Let us denote by $h(N)$ the left $A_{d+n}[t_0]$-submodule of 
$A_{d+n}[t_0]^r$ generated by $h(P_1),\dots, h(P_l)$. 
Let $\G = \{Q_1(t_0),\dots,Q_k(t_0)\}$ be a Gr\"obner basis of $h(N)$ 
with respect to $\prec_H$ 
consisting of $F[\mvec]$-homogeneous operators.  
Then $\G(1) := \{Q_1(1),\dots,Q_k(1)\}$ is an $F[\mvec]$-Gr\"obner basis of $N$.   
\end{proposition}

Thus we have an algorithm of computing an $F[\mvec]$-Gr\"obner 
basis for an arbitrary shift vector $\mvec \in \Z^r$.
We can prove the following in the same way 
as \cite[Proposition 3.11]{OakuBfunction}

\begin{proposition} \label{prop:Fdivision}
Let $N$ and $Q_j(t_0)$ be as in Proposition \ref{prop:F-Gbase} and put 
$\Nsc := \Dsc_X\otimes_{A_{d+n}}N \subset \Dsc_X^r$.  
Then for any germ $P$ of $\Nsc$ at $p \in Y$, there exist germs 
$U_j$ of $\Dsc_X$ at $p$ such that 
$P = U_1Q_1(1) + \cdots + U_kQ_k(1)$ and 
$\ord_Y[\mvec](U_jQ_j(1)) \leq \ord_Y[\mvec](P)$ for $j=1,...,k$.
\end{proposition}

If the leading exponent of $P \in A_{d+n}[t_0]^r$ is 
$\lexp_H(P) = (\lambda,\tilde\alpha,i) \in \N \times L \times \Set{r}$, 
we define the leading position $\lp_H(P)$ of $P$ by $i$. 
For $\tilde\alpha, \tilde\beta \in L$ and $i \in \{1,...,r\}$, we put
\begin{eqnarray*}
(\tilde\alpha,i) \vee (\tilde\beta,i) &:=& 
(\max\{\tilde\alpha_1,\tilde\beta_1\},...,
\max\{\tilde\alpha_{2d+2n},\tilde\beta_{2d+2n}\},i),\\
(\tilde\alpha,i) + (\tilde\beta,i) &:=&
(\tilde\alpha + \tilde\beta, i).
\end{eqnarray*}
Let $N$ and $Q_j(t_0)$ be as in Proposition \ref{prop:F-Gbase} 
and put $\Lambda := \{(i,j)\mid 1 \leq i <j \leq k,\,\, 
\lp_H(Q_i(t_0)) = \lp_H(Q_j(t_0)) \}$.  
For $(i,j) \in \Lambda$, let $S_{ij}(t_0), S_{ji}(t_0) \in A_{d+n}[t_0]$ 
be monomials such that 
\begin{eqnarray*}
\mbox{}&
\lexp_H(S_{ji}(t_0)Q_i(t_0)) = \lexp_H(S_{ij}(t_0)Q_j(t_0)) = 
\lexp_H(Q_i(t_0)) \vee \lexp_H(Q_j(t_0)),
\\
\mbox{}&
\lcoef_H(S_{ji}(t_0)Q_i(t_0)) = \lcoef_H(S_{ij}(t_0)Q_j(t_0)).
\end{eqnarray*}
Then by the Buchberger algorithm, there exist $F$-homogeneous 
$U_{ijl}(t_0) \in A_{d+n}[t_0]$  so that we have
$$
S_{ji}(t_0)Q_i(t_0) - S_{ij}(t_0)Q_j(t_0) 
= \sum_{l=1}^k U_{ijl}(t_0)Q_l(t_0) 
$$
and either $U_{ijl}(t_0) \neq 0$ or else
$$
\lexp_H(U_{ijl}(t_0)Q_l(t_0)) \prec_H \lexp_H(Q_i(t_0)) \vee \lexp(Q_j(t_0))
$$ 
for each $l = 1,...,k$.  

The proof of the following proposition is similar to that of 
\cite[Theorem 3.13]{OakuBfunction}:

\begin{proposition} \label{prop:syzygy}
In the same notation as in Proposition \ref{prop:Fdivision}, 
the left $A_{d+n}^r$-module 
$$
{\rm Syz}(Q_1(1),\dots,Q_k(1)) 
:= \{(U_1,...,U_k) \in A_{d+n}^r \mid U_1Q_1(1) + \cdots + U_kQ_k(1) = 0\}
$$
is generated by $\{V_{ij}(1) \mid (i,j) \in \Lambda\}$ with 
$$  
V_{ij}(t_0) := (0,\dots,\stackrel{(i)}{S_{ji}(t_0)},\dots,
\stackrel{(j)}{-S_{ij}(t_0)},\dots,0) - (U_{ij1}(t_0),\dots,U_{ijk}(t_0)).
$$
\end{proposition}

Now let us describe an algorithm for computing a free resolution of $M$
which is adapted to the filtration $F_Y[\mvec]$ 
(cf.\ Proposition \ref{prop:AtoD}).  
Let $N$ be a left $A_{d+n}$-submodule of $A_{d+n}^r$ such that 
$M = A_{d+n}^r/N$.  

First, starting with a given $\mvec \in \Z^r$, 
let $\{P_1,...,P_{r_1}\}$ be an $F[\mvec]$-Gr\"obner basis of 
$N$ constructed as in Proposition \ref{prop:F-Gbase}.  
Put
$$
\mvec_1 := (\ord_F[\mvec](P_1),\dots,\ord_F[\mvec](P_{r_1})).  
$$
and define $\psi_1 : A_{d+n}^{r_1} \longrightarrow A_{d+n}^{r}$ by 
$$
\psi_1(Q_1,...,Q_{r_1}) := Q_1P_1 + \cdots + Q_{r_1}P_{r_1}.
$$   
Then we get a set of generators of the kernel ${\rm Ker}\,\psi_1$ 
by using Proposition \ref{prop:syzygy}.  

By the same procedure as above with $N$, $r$, and $\mvec$  replaced by 
${\rm Ker}\,\psi_1$, $r_1$, and $\mvec_1$ respectively, 
we obtain a homomorphism $\psi_2: A_{d+n}^{r_2} \rightarrow A_{d+n}^{r_1}$ 
so that 
${\rm Im}\,\psi_2 = {\rm Ker\,\psi_1}$.  
In view of Propositions \ref{prop:Fdivision} and \ref{prop:syzygy}, 
the sequence
$$
F_Y^k[\mvec_2](A_{d+n}^{r_2}) \stackrel{\psi_2}{\longrightarrow}
F_Y^k[\mvec_1](A_{d+n}^{r_1}) \stackrel{\psi_1}{\longrightarrow}
F_Y^k[\mvec](A_{d+n}^r) 
$$
is exact for any $k \in \Z$ with $\mvec_2 \in \Z^{r_2}$ defined by
$$
\mvec_2 := (\ord_F[\mvec_1](\psi_2(1,\dots,0)),
\dots,\ord_F[\mvec_1](\psi_2(0,\dots,1))).  
$$
Proceeding in the same way, we can obtain a free resolution 
(\ref{eq:A-resolution}) 
which is adapted to the $F[\mvec]$-filtration for any given $l \in \N$.

\section{The $b$-function of a $D$-module}
\label{section:b-function}
\setcounter{equation}{0}

Let $\Msc$ be a left coherent $\Dsc_X$-module on $X$.  
We assume that a left $A_{d+n}$-submodule $N$ of $A_{d+n}^r$ is 
given explicitly so that $\Msc = \Dsc_X \otimes_{A_{d+n}}M$ holds 
with $M := A_{d+n}^r/N$.   
Set $\Nsc := \Dsc_X \otimes_{A_{d+n}} N \subset \Dsc_X^r$.  
We fix an arbitrary shift vector $\mvec= (m_1,\dots,m_r) \in \Z^r$ and put
\begin{eqnarray*}
\gr_Y^k[\mvec](\Nsc) &:=& F_Y^k[\mvec](\Nsc)/F_Y^{k-1}[\mvec](\Nsc), \\
\gr_Y^k[\mvec](\Msc) &:=& F_Y^k[\mvec](\Msc)/F_Y^{k-1}[\mvec](\Msc).
\end{eqnarray*}
They are left $\gr_Y^0(\Dsc_X)$-modules.
The $F[\mvec]$-filtration and the associated graded module are 
defined also for $A_{d+n}$-modules. 
Moreover we have
$$
\gr_Y^k[\mvec](\Msc) = \gr_Y^k(\Dsc_X)\gr_Y^0[\mvec](\Msc).
$$

We put $\vartheta := t_1\pd{t_1} + \cdots + t_d\pd{t_d}$.  
This is the unique vector field modulo $F_Y^{-1}(\Dsc_X)$ 
that operates on $\Jsc_Y/\Jsc_Y^2$ as identity.
Let $\theta$ be a commutative variable corresponding to $\vartheta$. 

\begin{definition} \label{def:b-function}
\rm 
The {\em $b$-function} (or the {\em indicial polynomial}) 
$b(\theta,p) \in K[\theta]$ of $\Msc$ along $Y$ 
with respect to the filtration $F_Y[\mvec]$ at $p \in Y$ 
is the monic polynomial $b(\theta,p) \in K[\theta]$ in $s$ of 
the least degree, if any, that satisfies
$$ 
b(\vartheta,p)\gr_Y^0[\mvec](\Msc)_p = 0. 
$$
If such $b(\theta,p)$ exists, $\Msc$ is called {\em specializable} along $X$ 
at $p$.  If $\Msc$ is not specializable at $p$, we put 
$b(\theta,p) = 0$.
The {\em global $b$-function} $b(\theta)$ of $\Msc$ along $Y$ is defined 
to be the least common multiple of $b(\theta,p)$ with $p$ running through 
$Y$. 
\end{definition}

It is known that the specializability does not depend on the shift vector 
$\mvec$ while the $b$-function can depend on it (cf.\ \cite{LS}). 
It is also known that if $\Msc$ is holonomic, then $\Msc$ is specializable 
at any $p \in X$ (\cite{KKasymptotic},\cite{KKholonomic3},
\cite{LaurentBfunction}).  

First, we reduce to the case $r=1$.
For each $i=1,\dots,r$, let $\pi_i$ be the projection of  
$\gr_Y^0[\mvec](\Dsc_X^r)$ to the $i$-th component and put 
$$
\gr_Y^0[\mvec](\Nsc)^{(i)} := \{ P \in \gr_Y^0[\mvec](\Nsc) \mid
\pi_j(P) = 0 \mbox{ for $j = i+1,\dots,r$}\}.
$$
Note that $\gr_Y^0[\mvec](\Nsc)^{(i)}/\gr_Y^0[\mvec](\Nsc)^{(i-1)}$ 
can be regarded as a left ideal of $\gr_Y^0(\Dsc_X)$ by the projection 
to the $i$-th component.    
Then we get the following lemma:

\begin{lemma}
Under the above notation, $b(\vartheta,p)$ is a generator of the ideal 
$$ 
\bigcap_{i=1}^r \left( K[\vartheta] \,\cap\, 
(\gr_Y^0[\mvec](\Nsc)^{(i)}/\gr_Y^0[\mvec](\Nsc)^{(i-1)})_p \right).
$$
\end{lemma}

Let us now assume that the order
$\prec$ satisfies
\begin{equation} \label{eq:graded-well-order}
\mbox{
$(\tilde\alpha,i) \prec (\tilde\alpha',j)$ if $i<j$ for 
$\tilde\alpha,\tilde\alpha' \in L$ and $i,j\in \{1,\dots,r\}$.  
}
\end{equation}
Let $\G$ be an $F[\mvec]$-Gr\"obner basis of $N$ with respect to $\prec_F$ 
defined by $\prec$ as in Section \ref{section:Gbase}.  
Then 
$
\hat\Isc := \gr_Y[\mvec](\Nsc)^{(i)}/\gr_Y[\mvec](\Nsc)^{(i-1)}
$ 
is generated by 
$$
\{\pi_i(\sigma_Y[\mvec](P)) \mid P \in \G,\,\,
\sigma_Y[\mvec](P) \in \gr_Y[\mvec](\Nsc)^{(i)}\}.
$$  

Our next task is to compute the intersection 
$\hat\Isc \cap \Dsc_Y[t_1\pd{1},\dots,t_d\pd{d}]$. 
For this purpose, we introduce commutative indeterminates 
$v = (v_1,\dots,v_d)$ and $w = (w_1,\dots,w_d)$, 
and work with the ring $A_{d+n}[v,w]$.  
For an element $P$ of $A_{d+n}$ of the form 
$$
P = \sum_{(\mu,\nu,\alpha,\beta)\in L} a_{\mu\nu\alpha\beta}
t^\mu\pd{t}
^\nu x^\alpha\pd{x}^\beta,
$$
its multi-homogenization ${\rm mh}(P) \in A_{d+n}[v]$ is defined by 
$$ 
{\rm mh}(P) := \sum_{(\mu,\nu,\alpha,\beta) \in L}
a_{\mu\nu\alpha\beta}v_1^{\nu_1-\mu_1-\kappa_1}\cdots 
v_d^{\nu_d-\mu_d-\kappa_d} t^\mu\pd{t}^\nu x^\alpha\pd{x}^\beta
$$
with $\kappa_j := \min\{\nu_j-\mu_j\mid a_{\mu\nu\alpha\beta} \neq 0\}$.  
Let $\prec_{mh}$ be an order on 
$\N^d \times \N^d \times L \ni (\rho,\sigma,\tilde\alpha)$ 
defined by 
\begin{eqnarray} \label{eq:mh-order}
& &
\mbox{
$(\rho,\sigma,\tilde\alpha) \prec_{mh} (\rho',\sigma',\tilde\alpha')$
if and only if $|\rho+\sigma| < |\rho'+\sigma'|$
} \nonumber \\
& & \mbox{
or else $|\rho+\sigma| = |\rho'+\sigma'|$, 
$\tilde\alpha < \tilde\alpha'$
}
\end{eqnarray}
with an arbitrary well-order $<$ on $L$ satisfying (\ref{eq:monomial-order}).
Fixing an $i \in \{1,\dots,d\}$, we assign weight $1$ to  
$w_i, \pd{t_i}$, weight $-1$ to $v_i, t_i$, 
and weight $0$ to all the other variables.  
An element of $A_{d+n}[v,w]$ is said to be multi-homogeneous if 
it is homogeneous with respect to the weight above for each 
$i=1,...,d$.  Thus ${\rm mh}(P)$ is multi-homogeneous for any 
$P \in A_{d+n}^r$.   
Put $S_\kappa := S_{1\kappa_1}\cdots S_{d\kappa_d}$
for $\kappa = (\kappa_1,\dots,\kappa_d)\in \Z^d$ with 
$S_{ij} = \pd{t_i}^j$ if $j \geq 0$ and 
$S_{ij} := t_i^{-j}$ otherwise.  
Let $s = (s_1,\dots,s_d)$ be commutative indeterminates.
Assume that $P \in A_{d+n}$ is multi-homogeneous. 
Then we have
$$
S_\kappa P = Q(t_1\pd{t_1},\dots,t_d\pd{t_d},x,\pd{x})
$$
with some $Q(s_1,\dots,s_d,x,\pd{x}) \in A_n[s_1,\dots,s_d]$
and $\kappa \in \Z^d$.
We put
$$
\psi(P)(s_1,\dots,s_d) := Q(s_1,\dots,s_d).
$$

\begin{proposition} \label{prop:subring}
Let $\hat\Isc$ be a left ideal of $\gr_Y(\Dsc_X)$.
Let $\G_0$ be a finite subset of $\gr_Y(A_{d+n})$ which generates 
$\hat\Isc$.  
Let $\G_1$ be a Gr\"obner basis with respect to $\prec_{mh}$ of the ideal of
$A_{d+n}[v,w]$ generated by 
$$ 
\{{\rm mh}(P) \mid P \in \G_0\} \cup \{ 1-v_iw_i \mid i =1,\dots,d\}.  
$$  
We may assume that $\G_1$ consists of multi-homogeneous elements 
since so does the input.
Then the left ideal 
$\hat\Isc \cap \Dsc_Y[t_1\pd{t_1},\dots,t_d\pd{t_d}]$
of $\Dsc_Y[t_1\pd{t_1},\dots,t_d\pd{t_d}]$ 
is generated by 
$$
\G_2 := \{\psi(P)(t_1\pd{t_1},\dots,t_d\pd{t_d}) 
\mid P \in \G_1 \cap A_{d+n}\}.
$$
\end{proposition}

\begin{proof}
Let $P$ be an element of $\G_1 \cap A_{d+n}$.  
Since $P$ is multi-homogeneous and free of $v,w$, there exists
$\kappa \in \Z^d$ so that 
$Q :=\psi(P)(t_1\pd{t_1},\dots,t_d\pd{t_d}) = S_\kappa P$.  
By definition, $P$ 
belongs to the ideal generated by ${\rm mh}(\G_0)$ and $1-v_iw_i$ 
($i=1,\dots,d$).  
Setting $v_i = w_i = 1$, we know that $Q$ belongs to $\hat\Isc$.

Conversely, let $P$ be an arbitrary germ of 
$\hat\Isc \cap \Dsc_Y[t_1\pd{t_1},\dots,t_d\pd{t_d}]$ 
at $p \in Y$.  
Multiplying $P$ by a polynomial in $x$ which does not vanish at 
$p$, we may assume $P \in A_n[t_1\pd{t_1},\dots,t_d\pd{t_d}]$.  
In view of the definition of the multi-homogenization and the fact that
${\rm mh}(P) = P$, there exists 
$\rho \in \N^d$ so that $v^\rho P$ belongs to the ideal generated by 
${\rm mh}(\G_0)$.  
This implies that $P$ belongs to the ideal generated by 
${\rm mh}(\G_0)$ and $\{ 1-v_iw_i \mid i = 1,\dots,d\}$ since
$$
P = (1-v^\rho w^\rho)P + v^\rho w^\rho P
$$ and $(1-v^\rho w^\rho)$ belongs to the ideal generated by
$\{ 1-v_iw_i \mid i = 1,\dots,d\}$.  
Set $\G_1 \cap A_{d+n} = \{P_1,\dots,P_k\}$.  
Then by the definition of $\prec_{mh}$ and $\G_1$, there exist
$Q_1,\dots,Q_k \in A_{d+n}$ so that
$$
P = Q_1P_1 + \cdots + Q_kP_k.
$$
Since $P_1,\dots,P_k$ are multi-homogeneous as well as $P$, we may assume
that such is also the case with $Q_1,\dots,Q_k$.  
Hence there exist $\kappa^{(1)},\dots,\kappa^{(k)} \in \Z^d$ so that
$$
P = \psi(P) = Q_1S_{\kappa^{(1)}}\psi(P_1) + \cdots + 
Q_kS_{\kappa^{(k)}}\psi(P_k).
$$
This completes the proof. 
\end{proof}

Let $\hat\Isc$ be as in Proposition \ref{prop:subring}.  
Now we have obtained a set of generators $\G_2$ of 
$\hat\Isc \cap \Dsc_Y[t_1\pd{t_1},\dots,t_d\pd{t_d}]$. 
We identify each $t_i\pd{t_i}$ with $s_i$.  
Then from $\G_2$, we can compute a set of generators $\G_3$ of the ideal  
$\hat\Isc \cap \Osc_Y[s]$ of $\Osc_Y[s]$
by eliminating $\pd{x}$ by means of Gr\"obner basis in the Weyl algebra
(see e.g.\ \cite{OakuChar} for details).  
Then it is easy to obtain a subset $\G_4$ of $K[x,\theta]$ 
which generates the sheaf of ideals
$$
\Jsc := \hat\Isc \cap \Osc_Y[\theta]
      = \hat\Isc \cap \Osc_Y[s] \cap \Osc_Y[\theta]
$$
with $\theta = s_1 + \cdots + s_d$ again by Gr\"obner basis 
in the polynomial ring.  
Let us denote by $J$ the ideal of $K[x,\theta]$ generated by $\G_4$.  
Then $b(\theta,p)$ is a generator of
$$
\Jsc_p \cap K[\theta] = (\Osc_Y[\theta])_pJ \cap K[\theta].
$$

Our final task is to compute the $b$-function at each point of $Y$ 
by using the input $\G_4$.  
This is achieved by primary decomposition.  
Let us state the method in a more general setting, where we replace
the variable $\theta$ by the variables $s = (s_1,\dots,s_d)$ 
for the sake of generality:
So let $J$ be an arbitrary ideal of $K[x,s]$ whose generators are given.   
For each point $p$ of $Y$, put 
$$
B(J,p) := (\Osc_Y[s])_pJ \cap K[s], 
$$
which is an ideal of $K[s]$.   
Let $J = Q_1 \cap \cdots \cap Q_l$ be a primary decomposition in $K[x,s]$.  
Then by the flatness of $(\Osc_Y[s])_p$ over $K[x,s]$ we have
$$
B(J,p) = B(Q_1,p) \cap \cdots \cap B(Q_l,p).
$$
Each ideal on the right hand side can be computed easily by the following:
 
\begin{lemma} \label{lemma:primary}
Let $Q$ be a primary ideal of $K[x,s]$ and put 
$$
{\bf V}_Y(Q) := \{ x \in Y=K^n \mid f(x) = 0 \mbox{ for any 
$f \in Q \cap K[x]$}\}.
$$ 
Then we have
$$
B(Q,p) = \left\{ \begin{array}{ll} Q\cap K[s] & \mbox{if $p\in {\bf V}_Y(Q)$}\\
                             K[s] & \mbox{if $p \in Y\setminus {\bf V}_Y(Q)$}.
\end{array}\right.
$$
\end{lemma}

\begin{proof}
First assume $p \not\in {\bf V}_Y(Q)$.  
Then there exists $a(x) \in K[x]$ such that $a(p) \neq 0$.  
This implies that $B(Q,p) = K[s]$.  
Next assume $p \in {\bf V}_Y(Q)$ and $b(s) \in B(Q,p)$.  
Then there exists $a(x) \in K[x]$ so that $a(x)b(s) \in Q$
and $a(p) \neq 0$.  
Suppose $b(s)$ does not belong to $Q$.  
Then we have $a(x)^j \in Q$ with some $j \in \N$ since $Q$ is primary.
This implies $a(p)=0$, which is a contradiction.
Thus we have $B(Q,p) \subset Q \cap K[s]$.  
The converse inclusion is obvious.    
\end{proof}

\begin{lemma} \label{lemma:lcm}
Let $J$ be an ideal of $K[x,s]$.  Then we have
$$ \bigcap_{p\in Y} B(J,p) = J \cap K[s].$$ 
\end{lemma}

\begin{proof}
Let $J = Q_1 \cap \cdots \cap Q_l$ be a primary decomposition.  
Then by the preceding lemma, we have
\begin{eqnarray*}
\bigcap_{p\in Y} B(J,p) &=& \bigcap_{p\in Y}\bigcap_{j=1}^l B(Q_j,p)\\
&=& \bigcap \{ Q_j \cap K[s] \mid Q_j \cap K[x] \neq K[x]\}\\
&=& J \cap K[s]
\end{eqnarray*}
since $Q_j \cap K[s] = K[s]$ if and only if $Q_j \cap K[x] = K[x]$.
\end{proof}

Returning back to the ideal $J$ of $K[s,\theta]$ generated by $\G_4$, 
we have only to apply Lemma \ref{lemma:primary} with $s$ replaced by 
$\theta$.  
This gives us an algebraic stratification of $Y$ so that $b(s,p)$ is 
constant on each stratum as a function of $p$.  
Moreover, Lemma \ref{lemma:lcm} tells us that the global $b$-function of 
$\Msc$ along $Y$ is simply a generator of $J \cap K[\theta]$. 
Thus the algorithm is summarized as follows:

\begin{algorithm} \label{algorithm:b-function}
\rm
(The $b$-function of $\Msc := \Dsc_X\otimes_{A_{d+n}}M$)

\noindent
Input: $M = A_{d+n}^r/N$ with an $A_{d+n}$-submodule $N$ of 
$A_{n+d}^r$, and $\mvec \in \Z^r$. 

\begin{enumerate}
\item 
Compute a Gr\"obner basis $\G$ of $N$ with respect to the order 
$\prec_F$ that is defined through (\ref{eq:F-order}) by using 
$\mvec$ and an order 
$\prec$ satisfying (\ref{eq:monomial-order}) and (\ref{eq:graded-well-order}).
\item
For $i = 1$ to $r$ do
\begin{enumerate}
\item
$\G_i := \{\pi_i(\sigma_Y[\mvec](P)) \mid P \in \G,\,\,
\pi_j(\sigma_Y[\mvec](P)) = 0 \,\, \mbox{for any $j > i$}\}$.
\item
Let $\G_{i1}$ be a Gr\"obner basis of the left ideal of $A_{d+n}[v,w]$ 
generated by
$$
\{{\rm mh}(P) \mid P \in \G_i\} \cup \{ 1-v_iw_i \mid i =1,\dots,d\}
$$
with respect to an order $\prec_{mh}$ satisfying (\ref{eq:mh-order}).
\item
$\G_{i2} := \{\psi(P) \in A_n[s] \mid P \in \G_{i1} \cap A_{d+n}\}$
\item
Compute $J_i := \langle \G_{i2} \rangle \cap K[x,\theta]$ 
first by eliminating $\pd{x}$, 
then eliminating $s'_2,\dots,s'_d$ after substitution
$\theta = s_1+\cdots+s_d$ and $s'_j = s_j$ for $j=2,\dots,d$;
here $\langle \G_{i2} \rangle$ denotes the left ideal of $A_n[s]$ generated
by $\G_{i2}$.
\end{enumerate}
\item The global $b$-function $b(\theta)$ of $\Msc$ is the generator of 
$\bigcap_{i=1}^r (J_i \cap K[\theta])$.
\item
For $i=1$ to $r$ do
\begin{enumerate}
\item
Compute a primary decomposition $J_i = \bigcap_{j=1}^{l_i}Q_{ij}$ 
in $K[x,\theta]$.
\item
Compute (generators of) $Q_{ij} \cap K[s]$ 
and $Q_{ij} \cap K[x]$ for $j = 1,\dots,l_i$ by elimination.
\end{enumerate}
\item For each $p \in Y$, the local $b$-function $b(\theta,p)$ is 
the generator of the ideal
$$ 
\bigcap_{i=1}^r \bigcap_{j=1}^{l_i}
\{ Q_{ij} \cap K[\theta] \mid 
g(p) = 0\,\, \mbox{for any $g(x) \in Q_{ij} \cap K[x]$}\}.
$$
\end{enumerate}
\end{algorithm}

Let us remark on the coefficient field:
Suppose that the input is defined over a subfield $K_0$ of $K$.  
Then the steps 1--3 can be done over $K_0$ instead of $K$ and 
$b(\theta,p)$ divides $b(\theta)$ for any $p \in Y$.   
However, the primary decomposition in the step 4 must be one in 
$K[x,\theta]$ not in $K_0[x,\theta]$.  
In fact, we need a primary decomposition over an intermediate field $K_1$ 
with $K_0 \subset K_1 \subset K$ so that $b(\theta)$ factors into 
linear polynomials in $K_1[\theta]$.  
If, e.g., $K_0$ is the rationals $\Q$, such $K_1$ is computable.  
Hence the primary decomposition in the step 4 is certainly computable 
if the input is defined over $\Q$ in view of 
e.g., \cite{BWBook}, \cite{EHV}, \cite{SY} and gives the local $b$-function 
at any $p \in Y = K^n$.

As a special case where all the computation can be done over $K_0$, 
suppose that the ideal $J \cap K_0[\theta]$ is generated by a polynomial 
which is a multiple of linear factors over $K_0$.
Then the step 4 of Algorithm \ref{algorithm:b-function} can be computed
over $K_0$ and the step 5 is true for any $p \in K^n$;
one can easily verify this by considering a projection of $K$ to $K_0$.  
Note that this is exactly the case with the classical Bernstein-Sato 
polynomial (cf.\ \cite{OakuBfunction},\cite{OakuMega96},\cite{OakuAdvance} 
for algorithms) and $K_0 = \Q$ 
by virtue of Kashiwara's theorem on the rationality 
\cite{KashiwaraBfunction}.  

At this occasion, let us make a correction to \cite{OakuAdvance}:
Lemma 4.4 of \cite{OakuAdvance} does not hold in general; 
we need field extension as explained above.
This correction does not affect the rest of \cite{OakuAdvance}.

\begin{example}\rm
Put $X := K^5 \ni (t_1,t_2,x,y,z)$ and 
$Y := \{(t_1,t_2,x,y,z) \in X \mid t_1 = t_2 = 0\}$.  
Put $M := A_5/I$ with the left ideal $I$ generated by
$$ 
t_1-x^3+y^2,\,\,
t_2-y^3+z^2,\,\,
\pd{x} + 3x^2\pd{t_1},\,\,
\pd{y} - 2y\pd{t_1} + 3y^2\pd{t_2},\,\,
\pd{z}-2z\pd{t_2}.
$$
Then the $b$-function $b(s,p)$ of $\Msc := \Dsc_X\otimes_{A_5}M$ along 
$Y$ at $p \in Y = K^3$ is given as follows:  
$$
\begin{array}{ll}
b(s,p) = &
s \left(s-\frac{5}{18}\right)\left(s-\frac{1}{6}\right)
\left(s-\frac{1}{18}\right)\left(s+\frac{1}{18}\right)
\left(s+\frac{1}{6}\right)\\
&\times
\left(s+\frac{5}{18}\right)
\left(s+\frac{1}{3}\right)\left(s+\frac{7}{18}\right)
\left(s+\frac{11}{18}\right)\left(s+\frac{2}{3}\right)
\end{array}
$$
if $p = (0,0,0) \in Y$; $b(s,p) = s$ if 
$p \in \{(x,y,z) \mid x^3-y^2 = 0,\,\, y^3-z^2 = 0 \}$; and 
$b(s,p) = 1$ otherwise. 
\end{example}

\section{Restriction of a $D$-module}
\label{section:restriction}
\setcounter{equation}{0}

We retain the notation of the preceding sections.  
Put 
$$
\Dsc_{Y\rightarrow X} := \Osc_Y \otimes_{\Osc_X}\Dsc_X.
$$ 
Then $\Dsc_{Y\rightarrow X}$ has a natural structure of 
$(\Dsc_Y, \Dsc_X)$-bimodule. 
Let $\Msc = \Dsc_X\otimes_{A_{d+n}}M$ be a coherent $\Dsc_X$-module 
with a finitely generated $A_{d+n}$-module $M = A_{d+n}^r/N$.  
Then the ($D$-module theoretic) {\em restriction} of $\Msc$ to $Y$ 
is defined by
$$
\Msc_Y^\bullet := \Dsc_{Y\rightarrow X}\otimes^L_{\Dsc_X}\Msc
$$
in the derived category of left $\Dsc_X$-modules 
(see \cite{HartshorneNotes} for the derived category and derived functors).  

In general, let $\Lsc$ be a $K[t]$-module and $\Lsc_j$ 
($j \in \Z$) be additive subgroups of $\Lsc$ such that 
$t_i\Lsc_j \subset \Lsc_{j-1}$ holds for $i = 1,\dots,d$ and $j \in \Z$.  
Then for any integer $k$, we define the Koszul complex associated with 
$\Lsc_\bullet = \{\Lsc_j\}_{j\in\Z}$ and $t_1,\dots,t_d$ by
$$ 
0 \longrightarrow 
\Lsc_{k+d} \otimes_\Z \stackrel{0}{\wedge}\Z^d 
\stackrel{\delta}{\longrightarrow}
\Lsc_{k+d-1}\otimes_\Z \stackrel{1}{\wedge}\Z^d 
\stackrel{\delta}{\longrightarrow}
\cdots
\stackrel{\delta}{\longrightarrow}
\Lsc_k \otimes_\Z \stackrel{d}{\wedge}\Z^d 
\longrightarrow 0,
$$
where $\delta$ is defined by 
$$
\delta(u\otimes e_{i_1}\wedge\cdots\wedge e_{i_j}) 
= \sum_{l=1}^d t_lu\otimes e_l\wedge e_{i_1}\wedge\cdots\wedge e_{i_j}
$$
for a subset $\{i_1,\dots,i_j\}$ of $\{1,\dots,d\}$ 
with the unit vectors $e_1,\dots,e_d$ of $\Z^d$.   
We denote this complex by $\Ksc^\bullet(\Lsc_\bullet[k],t_1,\dots,t_d)$.
When $\Lsc_j = \Lsc$ for each $j$, we also denote it simply by 
$\Ksc^\bullet(\Lsc,t_1,\dots,t_d)$. 
Here we regard $\Lsc_{k+d-j}\otimes\stackrel{j}{\wedge}\Z^d$ 
as being placed at the degree $-j$ to be compatible 
with the cohomology theory.

In particular, $\Ksc^\bullet(\Dsc_X,t_1,\dots,t_d)$ 
is quasi-isomorphic to $\Dsc_{Y\rightarrow X}$ 
in the derived category of right $\Dsc_X$-modules.
Hence we can identify $\Msc_Y^\bullet$ with the complex 
$$
\Ksc^\bullet(\Dsc_X,t_1,\dots,t_d)\otimes_{\Dsc_X}\Msc 
= \Ksc^\bullet(\Msc,t_1,\dots,t_d).  
$$

Our purpose below is to describe an algorithm to compute each cohomology 
group $\Hsc^i(\Msc_Y^\bullet)$ (for $i = 0,-1,\dots,-d$ since it is zero 
for other $i$) under the assumption that $\Msc$ is specializable along $Y$.
Let $b(\theta)$ be the global $b$-function of $\Msc$ with respect to 
the filtration $F_Y[\mvec]$ with a given $\mvec \in \Z^r$.  
In what follows, we can replace $b(\theta)$ by the local $b$-function 
$b(\theta,p)$ in order to compute $\Msc_Y^\bullet$ locally, i.e., 
on a Zariski neighborhood of $p$.  

\begin{proposition} \label{prop:criterion}
Let $k$ be an integer such that $b(k)\neq 0$.   
Then the Koszul complex \\
$\Ksc^\bullet(\gr_Y^\bullet[\mvec](\Msc)[k],t_1,\dots,t_d)$ 
associated with $\{\gr_Y^j[\mvec](\Msc)\}_{j\in\Z}$ 
is exact. 
\end{proposition}

We shall prove this proposition in a slightly more general situation.  
Let $A_d := K[t]\langle \pd{t} \rangle$ be the Weyl algebra 
on the variables $t = (t_1,\dots,t_d)$ and define a filtration on it and 
the associated graded module by
$$
F_k(A_d) := \{\sum_{\mu,\nu\in\N^d} a_{\mu\nu}t^\mu\pd{t}^\nu
\mid |\nu-\mu| \leq k\}, 
\qquad
\gr_k(A_d) := F_k(A_d)/F_{k-1}(A_d).
$$
Note that $\gr(A_d) := \bigoplus_{k\in\Z}\gr_k(A_d)$ is isomorphic to 
$A_d$.  In particular, we can regard $\gr_0(A_d)$ as a subring of $A_d$.

\begin{proposition} \label{prop:Koszul}
Let $\Lsc = \bigoplus_{j\in\Z}\Lsc_j$ be a graded $\gr(A_d)$-module; 
i.e., assume $\gr_j(A_d)\Lsc_i \subset \Lsc_{i+j}$ for $i,j \in\Z$.
Assume moreover that there exists a nonzero polynomial 
$b(\theta) \in K[\theta]$ which 
satisfies $b(\vartheta+j)\Lsc_j = 0$ for any $j\in\Z$ with 
$\vartheta = t_1\pd{t_1} + \cdots + t_d\pd{t_d}$.  
Let $k$ be an integer such that $b(k) \neq 0$.
Then $\Ksc(\Lsc_\bullet[k],t_1,\dots,t_d)$ is exact.  
\end{proposition}

\begin{proof}
We argue by induction on $d$.
First suppose $d=1$.  Then $\Ksc^\bullet(\Lsc_\bullet[k],t_1)$ is the complex
$$
0 \longrightarrow \Lsc_{k+1} \stackrel{t_1}{\longrightarrow}\Lsc_k 
\longrightarrow 0.
$$
Assume $u \in \Lsc_{k+1}$ satisfies $t_1u= 0$.  Then we have $u=0$ since
$$
0 = b(t_1\pd{t_1}+k+1)u = b(\pd{t_1}t_1 + k)u = b(k)u.
$$ 
On the other hand, there exists $P \in A_d$ so that 
$b(t_1\pd{t_1}+k) = t_1P + b(k)$.  
Hence for an arbitrary $v \in\Lsc_k$,  we conclude
$v \in t_1\Lsc_{k+1}$ from $b(t_1\pd{t_1} + k)v = 0$.

Now assume the proposition is true with $d$ replaced by $d-1$.  
It is easy to see, as in the case of the usual Koszul complex
(see e.g., \cite[p.188]{Schapira}), that 
$\Ksc^\bullet(\Lsc_\bullet[k],t_1,\dots,t_d)$ is quasi-isomorphic to 
the complex associated with the double complex
\begin{equation} \label{eq:doubleComplex}
\begin{array}{c}
\Ksc^\bullet(\Lsc_\bullet[k+1],t_1,\dots,t_{d-1})\\
\downarrow \, t_d \\
\Ksc^\bullet(\Lsc_\bullet[k],t_1,\dots,t_{d-1}).
\end{array}
\end{equation}
Let us denote by $\Lsc'_j$ and $\Lsc''_j$ the kernel and the cokernel of
$t_d : \Lsc_{j+1} \longrightarrow \Lsc_j$. 
Then $\Lsc' := \bigoplus_{j\in\Z}\Lsc'_j$ and 
$\Lsc'' := \bigoplus_{j\in\Z}\Lsc''_j$ are graded $\gr(A_{d-1})$-modules.  
For any $u \in \Lsc'_j$, we have
\begin{eqnarray*}
0 &=& b(t_1\pd{t_1} + \cdots + t_d\pd{t_d}+j+1)u \\
  &=& b(t_1\pd{t_1} + \cdots + t_{d-1}\pd{t_{d-1}}
       +\pd{t_d}t_d +j)u
\\
 &=& b(t_1\pd{t_1} + \cdots + t_{d-1}\pd{t_{d-1}}+j)u.
\end{eqnarray*}
On the other hand, for $v \in \Lsc_j$, let $\overline v$ be its 
residue class in $\Lsc''_j$.  Then we have
$$
0 = b(t_1\pd{t_1} + \cdots + t_d\pd{t_d}+j)\overline v
  = b(t_1\pd{t_1} + \cdots + t_{d-1}\pd{t_{d-1}}+j)\overline v.
$$
Thus both $\Lsc'$ and $\Lsc''$ satisfy the conditions of the proposition 
with $d$ replaced by $d-1$.  
By the induction hypothesis, the complexes 
$\Ksc^\bullet(\Lsc'_\bullet[k];t_1,\dots,t_d)$ and
$\Ksc^\bullet(\Lsc''_\bullet[k];t_1,\dots,t_d)$ 
are exact.
Hence the vertical chain map of 
(\ref{eq:doubleComplex}) is a quasi-isomorphism,  
which implies that $\Ksc^\bullet(\Lsc_\bullet[k],t_1,\dots,t_d)$ is exact. 
\end{proof}

Under the assumption of Proposition \ref{prop:criterion}, we have
$b(\vartheta+j)\gr_Y^j[\mvec](\Msc) = 0$ for any $j\in\Z$.  
In fact, for $P \in \gr_Y^j(\Dsc_X)$, we easily get
$b(\vartheta+j)P = Pb(\vartheta)$.  This yields
\begin{eqnarray*}
b(\vartheta+j)\gr_Y^j[\mvec](\Msc)
&=& b(\vartheta+j)\gr_Y^j(\Dsc_X)\gr_Y^0[\mvec](\Msc)\\
&=& \gr_Y^j(\Dsc_X)b(\vartheta)\gr_Y^0[\mvec](\Msc) = 0.
\end{eqnarray*}
Hence Proposition \ref{prop:criterion} is an immediate consequence of 
Proposition \ref{prop:Koszul}.

Now for $\mvec \in \Z^r$, 
we define the $F_Y[\mvec]$-filtration on $\Dsc_{Y\rightarrow X}^r$ by
\begin{eqnarray*}
F_Y^k[\mvec](\Dsc_{Y\rightarrow X}^r) 
&:=&
F_Y^k[\mvec](\Dsc_X^r)/
(t_1F_Y^{k+1}[\mvec](\Dsc_X^r) + \cdots + t_d F_Y^{k+1}[\mvec](\Dsc_X^r))
\\ 
&\simeq& 
\{ P = \sum_{i=1}^r\sum_{\nu,\beta} 
a_{\nu\beta}(x)\pd{t}^\nu\pd{x}^\beta e_i 
\mid a_{\nu\beta}(x) = 0\,\,\mbox{if $|\nu| > k-m_i$} \}\\
&=&
\bigoplus_{i=1}^r \bigoplus_{|\nu| \leq k-m_i} \Dsc_Y.  
\end{eqnarray*}

\begin{theorem} \label{th:truncation}
For an arbitrary $\mvec \in \Z^r$,  
let us take a free resolution (\ref{eq:resolution}) of Definition 2.1 
of length $l = d+1$, which is adapted to the $F_Y[\mvec]$-filtration.  
Let $\mvec_1,\dots,\mvec_{d+1}$ be the associated shift vectors.  
Take two integers $k_0 \leq k_1$ 
so that the global $b$-function $b(\theta)$ of $\Msc$ satisfies 
$b(k)\neq 0$ for any integer $k$ with $k < k_0$ or else $k > k_1$.
Then $\Msc_Y^\bullet$ is quasi-isomorphic to the complex
\begin{equation} \label{eq:truncation}
\cdots \rightarrow
\frac{F_Y^{k_1}[\mvec_{d+1}](\Dsc_{Y\rightarrow X}^{r_{d+1}})}
{F_Y^{k_0-1}[\mvec_{d+1}](\Dsc_{Y\rightarrow X}^{r_{d+1}})}
\stackrel{\overline\psi_{d+1}}{\longrightarrow}
\cdots
\stackrel{\overline\psi_1}{\longrightarrow}
\frac{F_Y^{k_1}[\mvec_0](\Dsc_{Y\rightarrow X}^{r_0})}
{F_Y^{k_0-1}[\mvec_0](\Dsc_{Y\rightarrow X}^{r_0})}
\rightarrow 0
\end{equation}
with $r_0 = r$ and $\mvec_0 = \mvec$, where $\overline\psi_j$ is 
a homomorphism induced by $\psi_j$.   
In particular, we have $\Msc_Y^\bullet = 0$ if $b(k) \neq 0$ 
for any $k\in\Z$.
\end{theorem}

\begin{proof}
For any $k \in \Z$, the complex 
$$
\cdots \rightarrow
F_Y^k[\mvec_{d+1}](\Dsc_X^{r_{d+1}}) 
\stackrel{\psi_{d+1}}{\longrightarrow} 
\cdots 
\stackrel{\psi_1}{\longrightarrow} F_Y^k[\mvec_0](\Dsc_X^{r_0})
\rightarrow 0
$$
is quasi-isomorphic to $F_Y^k[\mvec](\Msc)$ in view of the 
exact sequence (\ref{eq:resolution}) of Definition \ref{def:resolution}.  
Hence we know that 
$\Ksc^\bullet(\gr_Y^\bullet[\mvec](\Msc)[k],t_1,\dots,t_d)$
is quasi-isomorphic to the complex associated with the double complex
$$
\begin{array}{c}
\def\upmap#1{\rule[-12pt]{0pt}{30pt}\Big\uparrow\rlap
  {$\vcenter{\hbox{$\scriptstyle#1$}}$}}
0 \\
\upmap{} \\
\Ksc^\bullet(\gr_Y^\bullet[\mvec_0](\Dsc_X^{r_0})[k],t_1,\dots,t_d) \\
\upmap{\overline\psi_1} \\
\vdots \\
\upmap{\overline\psi_{d+1}}\\
\Ksc^\bullet(\gr_Y^\bullet[\mvec_{d+1}](\Dsc_X^{r_{d+1}})[k],t_1,\dots,t_d)\\
\upmap{} \\
\vdots
\end{array}
$$
On the other hand, we have a quasi-isomorphism
$$ 
\Ksc^\bullet(\gr_Y^\bullet[\mvec_i](\Dsc_X^{r_i})[k],t_1,\dots,t_d)
\,\,\simeq\,\,
\gr_Y^k[\mvec_i](\Dsc_{Y\rightarrow X}^{r_i}).
$$  
Hence $\Ksc^\bullet(\gr_Y^\bullet[\mvec](\Msc)[k],t_1,\dots,t_d)$
is quasi-isomorphic to the complex
\begin{equation} \label{eq:grComplex}
\cdots \rightarrow
\gr_Y^{k}[\mvec_{d+1}](\Dsc_{Y\rightarrow X}^{r_{d+1}})
\stackrel{\overline\psi_{d+1}}{\longrightarrow}
\cdots
\stackrel{\overline\psi_1}{\longrightarrow}
\gr_Y^{k_1}[\mvec_0](\Dsc_{Y\rightarrow X}^{r_0})
\rightarrow 0.   
\end{equation}
By virtue of Proposition \ref{prop:criterion}, 
the complex (\ref{eq:grComplex}) is exact if $b(k) \neq 0$.  
This implies the theorem since we have 
$$
F_Y^k[\mvec_i](\Dsc_{Y\rightarrow X}^{r_i}) = 0
$$
for sufficiently small $k\in\Z$.
\end{proof}

Note that 
$
F_Y^{k_1}[\mvec_i](\Dsc_{Y\rightarrow X}^{r_i})/
F_Y^{k_0-1}[\mvec_i](\Dsc_{Y\rightarrow X}^{r_i})
$
is a free $\Dsc_Y$-module of rank 
$$
\sum_{j=1}^{r_i}
\sharp\{\nu \in \Z^d \mid k_0-m_{ij} \leq |\nu| \leq k_1-m_{ij}\}.
$$
Hence Theorem \ref{th:truncation}
 gives us an algorithm to compute each cohomology group 
$\Hsc^i(\Msc_Y^\bullet)$.  In fact, we have only to compute the cohomology 
groups of the complex (\ref{eq:truncation}) as left $A_n$-modules 
with $\Dsc_{Y\rightarrow X}$ replaced by $A_n[\pd{t}]$.  
The flatness of $\Dsc_Y$ over $A_n$ assures us that the generators of 
the cohomology group over $A_n$ also give the ones over $\Dsc_Y$.  
The algorithm is summarized as follows: 

\begin{algorithm} \label{algorithm:restriction}
\rm
(The cohomology groups of the restriction of $\Msc$ to $Y$)

\noindent
Input: $M = A_{d+n}^r/N$ with an $A_{d+n}$-submodule $N$ of 
$A_{d+n}^r$. \\
Output: $\Hsc^{-i}(\Msc_Y^\bullet) = 
\Dsc_Y\otimes_{A_n}(A_n^{l_i}/I_i)$ for $0 \leq i \leq d$.  

\begin{enumerate}
\item
Choose an arbitrary $\mvec \in \Z^r$; e.g., take $\mvec = (0,\dots,0)$ 
by default.  
\item
Compute the global $b$-function $b(\theta)$ of $\Msc$ along $Y$ 
by the steps 1--3 of Algorithm \ref{algorithm:b-function} 
with $M$ and $\mvec$ as input.
\item If $b(\theta) = 0$, then $\Msc$ is not globally specializable along $Y$;
quit.
\item Let $k_0$ and $k_1$ be the minimum and the maximum integral root of 
$b(\theta) = 0$.  If there is no integral root, then we have 
$\Hsc^i(\Msc_Y^\bullet) = 0$ for all $i$; quit. 
\item
Compute a free resolution 
$$
A_{d+n}^{r_{d+1}} \stackrel{\psi_{d+1}}{\longrightarrow} 
A_{d+n}^{r_{d}} 
\stackrel{\psi_{d}}{\longrightarrow} 
\cdots 
\stackrel{\psi_2}{\longrightarrow} A_{d+n}^{r_1}
\stackrel{\psi_1}{\longrightarrow} A_{d+n}^{r_0}
\stackrel{\varphi}{\longrightarrow} M \longrightarrow 0
$$
of $M$ adapted to the $F_Y[\mvec]$-filtration 
(cf.\ Proposition \ref{prop:AtoD})
and the shift vectors $\mvec_1,\dots,\mvec_{d+1}$ 
by using Proposition \ref{prop:syzygy} successively, or by using 
Theorem \ref{th:Schreyer}.
\item 
Compute the induced complex
$$
\frac{F_Y^{k_1}[\mvec_{d+1}](A_n[\pd{t}]^{r_{d+1}})}
{F_Y^{k_0-1}[\mvec_{d+1}](A_n[\pd{t}]^{r_{d+1}})}
\stackrel{\overline\psi_{d+1}}{\longrightarrow}
\cdots
\stackrel{\overline\psi_1}{\longrightarrow}
\frac{F_Y^{k_1}[\mvec_0](A_n[\pd{t}]^{r_0})}
{F_Y^{k_0-1}[\mvec_0](A_n[\pd{t}]^{r_0})}
\rightarrow 0
$$
as a complex of finitely generated free left $A_n$-modules,
where $A_n[\pd{t}]$ is identified with 
$A_{d+n}/(t_1A_{d+n}+\cdots+ t_dA_{d+n})$. 
Put $\overline \psi_0 := 0$.  
\item 
Via Gr\"obner bases of modules over $A_n$,
compute the $-i$-th cohomology group
${\rm Ker}\,\overline \psi_i/{\rm Im}\,\overline \psi_{i+1}$
of the above complex in the form $A_n^{l_i}/I_i$ 
with a left $A_n$-module $I_i$ for $i=0,\dots,d$. 
\end{enumerate}
\end{algorithm}

Note that in the step 5 of the above algorithm, 
only $\psi_1,\dots,\psi_{i_0+1}$ 
are needed if one wants to compute only the $-i$-th cohomology groups 
for $i=0,\dots,i_0$.  In particular, one does not need the free resolution 
to compute only the $0$-th cohomology as will be the case with Algorithm 
\ref{algorithm:localization}. 

As a direct application of the algorithm above, we obtain an algorithm to
compute the cohomology groups with coefficients in the formal power 
series solutions of $\Msc$: 
$$
\Ext^i_{\Dsc_X}(\Msc,K[[x]]) \quad (i=0,\dots,n)
$$
under the assumption that $\Msc$ is specializable along $Y := \{0\}$.
In fact, we can easily verify that there exists an isomorphism 
(see e.g.\ \cite[p.428]{LM} for the case $K=\C$)  
$$
\Ext^i_{\Dsc_X}(\Msc,K[[x]]) \simeq \Ext^i_K(\Msc_Y^\bullet,K)
\simeq H^i(\Msc_Y^\bullet).
$$
If $K = \C$ and $\Msc$ is Fuchsian along $Y$ in the sense of \cite{LM}
(this condition holds if $\Msc$ is regular holonomic in the sense of 
\cite{KKholonomic3}), 
then by virtue of the comparison theorem, we have also an isomorphism
$$
\Ext^i_{\Dsc_X}(\Msc,\C\{x\}) \simeq H^i(\Msc_Y^\bullet),
$$
where $\C\{x\}$ denotes the ring of convergent power series in $x$.
Note that Kashiwara's index theorem (\cite[p.127]{KashiwaraBook})
gives the local index
$$
\sum_{i\geq 0}(-1)^i \dim_{\C}\Ext^i_{\Dsc_X}(\Msc,\C\{x\})
$$
at $0 \in X$ in terms of some topological quantity associated with the
characteristic cycle of $\Msc$.

\begin{example}\rm
Let us consider $M := A_4/I$, where $I$ is the left ideal of $A_4$ 
(with $K = \C$) generated by 
$$
x_3\pd{3} + x_4\pd{4} - a_1, \quad
x_1\pd{1} + x_3\pd{3} - a_2,\quad
x_2\pd{2} + x_4\pd{4} -a_3,\quad
\pd{1}\pd{4} - \pd{2}\pd{3}, 
$$
where $a_1,a_2,a_3 \in \C$ are parameters.  
Put $\Msc := \Dsc_X\otimes_{A_4}M$ with $X = \C^4$ and 
$Y_1 := \{(x_1,x_2,x_3,x_4) \in X \mid x_1 = x_2 = x_3 = 0\}$.
The global $b$-function of $\Msc$ along $Y_1$ is 
$(s-a_2)(s+a_1-a_2-a_3)$.  
Hence the cohomology groups of the restriciton of $\Msc$ to $Y_1$ all 
vanish unless $a_2$ or $a_1-a_2-a_3$ is an integer.  
If $a_1 = a_2 = a_3 = 0$, we have by Algorithm \ref{algorithm:restriction} 
$$
\Hsc^i(\Msc^\bullet_{Y_1}) = \left\{\begin{array}{ll}
\Dsc_{Y_1}/\Dsc_{Y_1}x_4\pd{4}          & (i=0),\\
(\Dsc_{Y_1}/\Dsc_{Y_1} x_4\pd{4})^2      & (i=-1),\\
\Dsc_{Y_1}/\Dsc_{Y_1} x_4\pd{4}          & (i=-2),\\
0                                & (i \leq -3).
\end{array}\right.
$$
The $b$-function of $\Msc$ along the point $Y_0 := \{(0,0,0,0)\}$ is 
$s-a_2-a_3$.  
Suppose $a_1 = a_2 = a_3 = 0$.  
Then the cohomology groups $\Hsc^i(\Msc^\bullet_{Y_0})$ 
of the restriction of $\Msc$ to $Y_0$ are
$\C,\C^3,\C^3,\C,0$ for $i = 0,-1,-2,-3,-4$ respectively.  
Since $\Msc$ is regular holonomic, this implies that 
$\Ext_{\Dsc_X}^i(\Msc,\C\{x\})$ is $\C,\C^3,\C^3,\C,0$ for $i=0,1,2,3,4$ 
respectively. 
\end{example}

\begin{example}\rm
Put $X := \C^2 \ni (x,y)$ and 
$\Msc := \Dsc_X\otimes_{A_2}(A_2/I)$ with $I$ being the left ideal generated by
$$
x\pd{x} - x(x\pd{x}+y\pd{y}+a)(x\pd{x}+b_1),\quad
y\pd{y} - y(x\pd{x}+y\pd{y}+a)(y\pd{y}+b_2).
$$
Then by the computation of the restriction of $\Msc$ to $(0,0)$, we get
$$
\Extsc^i_{\Dsc_X}(\Msc,\C[[x,y]]) = \left\{ \begin{array}{ll}
\C   & (i=0),\\
\C^2 & (i=1),\\
\C    & (i=2).
\end{array} \right.
$$
for generic parameters $a,b_1,b_2$ (this means that we perform the 
computation over the coefficient field $K := \Q(a,b_1,b_2)$).   
In particular, we have 
$\sum_{i=0}^2 (-1)^i\dim_{\C} \Extsc_{\Dsc_X}^i(\Msc,\C[[x,y]]) = 0$.  
On the other hand, the characteristic cycle of $\Msc$ is 
$$
3\{\xi=\eta=0\} + 4\{x=\eta = 0\} + 4\{y=\xi=0\} + \{x-y=\xi+\eta=0\} 
+ 7\{x = y = 0\}
$$
as a cycle in the cotangent bundle $T^*X = \{(x,y,\xi,\eta)\}$.  
Thus, by Kashiwara's index theorem we have
$$
\sum_{i=0}^2 (-1)^i\dim_{\C}\Extsc_{\Dsc_X}^i(\Msc,\C\{x,y\}) 
= 3 - (4+4+1) + 7 = 1.
$$
Hence $\Msc$ is not regular at $(0,0)$.   

${\rm dim}_{\C}\, \Extsc_{\Dsc_X}^0(\Msc,\C[[x,y]]) = 1$
implies that the system $\Msc$ admits one dimensional space of formal
power series solutions at the origin.
In fact, the (divergent) formal series
$$ \sum_{m,n = 0}^\infty 
   \frac{(a)_{m+n}(b_1)_m (b_2)_n}
        {(1)_m (1)_n} x^m y^n, \quad (c)_m:=c (c+1) \cdots (c+m-1)
$$
spans the solution space.
\end{example}

\section{Tensor product and localization}
\label{section:tensor}
\setcounter{equation}{0}

In this and subsequent sections, we denote by $X$ the affine space 
$K^n$.
First let us describe an algorithm to compute the tensor product and the
torsion groups of two holonomic $\Dsc_X$-modules $\Msc_1$ and $\Msc_2$. 
We suppose that left $A_n$-modules $N_1$ and $N_2$ are given so that
$$
M_i := \Gamma(X,\Msc_i) = A_n^{r_i}/N_i \quad (i=1,2).
$$
Let $\pi_1,\pi_2 : X \times X \rightarrow X$ be the projections to the
first and the second component respectively and put
$$
\Dsc'_{X\times X} := \pi_1^{-1}\Dsc_X\otimes_K\pi_2^{-1}\Dsc_X.  
$$
Then the exterior tensor product is defined by
$$
\Msc_1 \hat\otimes \Msc_2 := 
\Dsc_{X\times X}\otimes_{\Dsc'_{X\times X}}
(\pi_1^{-1}\Msc_1 \otimes_K \pi_2^{-1}\Msc_2).
$$

First, let us describe this exterior tensor product more concretely.
Let  $u_1,\dots,u_{r_1}$ be the residue classes of the unit
vectors $e_1,\dots,e_{r_1}$ of $A_n^{r_1}$, 
and $v_1,\dots,v_{r_2}$ the residue classes of the unit vectors 
$e'_1,\dots,e'_{r_2}$ of $A_n^{r_2}$.  
Then as a $\Dsc_{X\times X}$-module, 
$\Msc_1\hat\otimes \Msc_2$ is generated by $u_i\otimes v_j$ 
($1 \leq i \leq r_1$, $1 \leq j \leq r_2$).
Let us denote by $(x,y)$ the coordinate system of $X\times X$. 
For $P = (P_1,\dots,P_{r_1}) \in \Dsc_X^{r_1}$ 
and $Q = (Q_1,\dots,Q_{r_2}) \in \Dsc_X^{r_2}$ 
we write
$$
P\otimes Q := (P_i(x,\pd{x})Q_j(y,\pd{y}))_{ij} 
\in\Dsc_{X\times X}^{r_1r_2}.
$$
Let $\Isc'$ be the left $\Dsc'_{X\times X}$-submodule
of $(\Dsc'_{X\times X})^{r_1r_2}$ generated by the set
$$
\{ P\otimes e'_j,\,\, e_i\otimes Q \mid  
P \in N_1,\,\, Q \in N_2,\,\, 
1 \leq i \leq r_1, \,\, 1 \leq j \leq r_2
\}
$$
and put
$\Isc := \Dsc_{X\times X}\otimes_{\Dsc'_{X\times X}}\Isc'$.

\begin{lemma} \label{lemma:exterior-tensor-product}
Under the above notation, there is an isomorphism
$$
\Msc_1\hat\otimes \Msc_2 \,\,\simeq\,\,
\Dsc_{X\times X}^{r_1r_2}/\Isc.
$$ 
\end{lemma} 

\begin{proof}
Put 
$$
\Ksc' := \{(P_{ij})_{ij} \in (\Dsc'_{X\times X})^{r_1r_2}
\mid \sum_{i,j} P_{ij}(u_i\otimes v_j) = 0\,\,
\mbox{in $\pi_1^{-1}\Msc_1\otimes_K\pi_2^{-1}\Msc_2$} \}.  
$$ 
Then it is easy to see that $\Isc' \subset \Ksc'$.
Hence we have a commutative diagram
$$
\matrix{ \pi_1^{-1}\Msc_1 \times \pi_2^{-1}\Msc_2 & 
\stackrel{\Phi}\longrightarrow & 
(\Dsc'_{X\times X})^{r_1r_2}/\Isc' & \cr
       &  \searrow & \downarrow & \cr
& & \pi_1^{-1}\Msc_1 \otimes_K \pi_2^{-1}\Msc_2 & 
\simeq (\Dsc'_{X\times X})^{r_1r_2}/\Ksc', \cr
}
$$
where $\Phi$ is a $K$-bilinear map defined by
$$
\Phi(\sum_i P_iu_i, \sum_jQ_jv_j) = P\otimes Q \mod \Isc',
$$
which is well-defined by the definition of $\Isc'$.
In view of the universal property of the tensor product, we know 
that the vertical map above is an isomorphism.
This completes the proof.
\end{proof}

Hence $\Msc_1\hat\otimes\Msc_2$ is computable with $N_1$ and 
$N_2$ being given.  Put
$\Delta := \{(x,y) \in X \times X \mid x=y\}$ and identify $\Delta$ 
and $X$ by the map $\pi_1$. 
Then by \cite[Proposition 4.7]{HolonomicII}, which obviously applies 
to algebraic $\Dsc$-modules as well, we have
\begin{eqnarray*}
\Msc_1 \otimes^L_{\Osc_X} \Msc_2 &\simeq& 
\Dsc_{\Delta\rightarrow X\times X}\otimes^L_{\Dsc_{X\times X}} 
(\Msc_1\hat\otimes\Msc_2) \\
&=&
(\Msc_1\hat\otimes\Msc_2)_\Delta^\bullet.  
\end{eqnarray*}
Suppose that $\Msc_1$ and $\Msc_2$ are holonomic.
Then it is easy to see that $\Msc_1\hat\otimes\Msc_2$ is a holonomic 
$\Dsc_{X\times X}$-module since its characteristic variety is contained in 
the Cartesian product of those of $\Msc_1$ and $\Msc_2$.  
Hence $\Msc_1\hat\otimes\Msc_2$ is specializable along $\Delta$ and 
the following algorithm is correct:

\begin{algorithm} \label{algorithm:tensor}
\rm
(The tensor product and torsion groups of two $\Dsc_X$-modules)

\noindent
Input: Holonomic systems $\Msc_i = \Dsc_X \otimes_{A_n}M_i$, where
$M_i = A_n^{r_i}/N_i$ with an $A_n$-submodule $N_i$ of 
$A_n^{r_i}$ for $i=1,2$. \\
Output: $\Torsc_k^{\Osc_X}(\Msc_1,\Msc_2) = \Lsc_k$ for $k=0,\dots,n$.  

\begin{enumerate}
\item From sets $\G_i$ of generators of $N_i$, compute
$$ 
\G_3 := \{ P\otimes e'_j,\,\, e_i\otimes Q \mid  
P \in \G_1,\,\, Q \in \G_2,\,\, 
1 \leq i \leq r_1, \,\, 1 \leq j \leq r_2
\}.
$$
\item 
Let $\G_4$ be the result of the substitution $y_i = x_i + t_i$ 
($i=1,\dots,n$) for each element of $\G_3$; 
let $I_4$ be the left ideal of 
$A_{2n} = K[t,x]\langle \pd{t},\pd{x}\rangle$ generated by $\G_4$ 
with $t = (t_1,\dots,t_n)$
\item
Apply Algorithm \ref{algorithm:restriction} with $A_{2n}/I_4$ 
as input and $d=n$ to obtain 
$$
\Lsc_k = \Hsc^{-k}((\Dsc_{X\times X}\otimes_{A_{2n}}
(A_{2n}^{r_1r_2}/I_4))_{\{0\}\times X}^\bullet).
$$
for $k = 0,\dots,n$.
\end{enumerate}
\end{algorithm}

\begin{example}\rm
Put $X := K$ and 
$$
\Msc := \Dsc_X/\Dsc_X x\pd{x},\quad \Nsc := \Dsc_X/\Dsc_X x.
$$
First, the exterior tensor product is given by
$$
\Msc\hat\otimes\Nsc =  \Dsc_{X\times X}/(
\Dsc_{X\times X}x\pd{x} + \Dsc_{X\times X}y) 
$$
with $(x,y) \in X\times X$.  Its global $b$-function along the diagonal 
is $s$, and by restricting $\Msc\hat\otimes\Nsc$ to the diagonal we get
$$
\Torsc^{\Osc_X}_i(\Msc,\Nsc) = \Dsc_X/\Dsc_X x \quad (i=0,1).
$$
In the same way, we get
\begin{eqnarray*}
\Torsc^{\Osc_X}_i(\Dsc_X/\Dsc_X x,\Dsc_X/\Dsc_X x) 
&=& \left\{\begin{array}{ll} 0 & (i=0),\\
\Dsc_X/\Dsc_X x & (i=1).\end{array}\right.\\
\Torsc^{\Osc_X}_i(\Dsc_X/\Dsc_X(x\pd{x}+1),\Dsc_X/\Dsc_X x) 
&=& 0 \quad (i=0,1).
\end{eqnarray*}
\end{example}

Let $\Msc$ be a holonomic $\Dsc_X$-module and let $f \in K[x]$ 
be an arbitrary non-constant polynomial.  
Then we immediately obtain an algorithm to compute the localization 
$\Msc[f^{-1}] := \Osc_X[f^{-1}]\otimes_{\Osc_X}\Msc$
by combining this algorithm with that of computing $\Osc[1/f]$ 
given in \cite{OakuMega96} since $\Osc_X[f^{-1}]$ is holonomic 
(\cite[Theorem 1.3]{HolonomicII}).  
Since $\Osc_X[f^{-1}]$ is flat over $\Osc_X$, the higher torsion 
groups vanish.  

\begin{algorithm} \label{algorithm:localization}
\rm
(The localization $\Msc[f^{-1}]$)

\noindent
Input: A holonomic system $\Msc = \Dsc_X \otimes_{A_n}(A_n^r/N)$
and a non-constant polynomial $f \in K[x]$. \\
Output: $\Msc[f^{-1}]$.

\begin{enumerate}
\item
Compute the global Bernstein-Sato polynomial $b_f(s)$ of $f$ as 
follows (cf.\ \cite{OakuBfunction}), where $s$ is a single indeterminate:
\begin{enumerate}
\item
Letting $t$ be a single variable and let $I$ be the left ideal of 
$A_{n+1}$ generated by 
$t-f(x)$ and $\pd{x_i} + (\partial f/\partial x_i)\pd{t}$ 
for $i=1,\dots,n$.
\item
Let $b(\theta)$ be that in the step 3 of Algorithm \ref{algorithm:b-function} 
with $A_{n+1}/I$ as input and $d=1$.
Put $b_f(s) := b(-s-1)$.
\end{enumerate}
\item
Compute a set of generators of 
$J_f := \{P(s,x,\pd{x}) \mid A_n[s] \mid Pf^s = 0\}$ 
by \cite[Theorem 19]{OakuMega96}. 
\item
Let $\nu$ be the minimum integer root of $b_f(s) = 0$ and put
$J_f(\nu) := \{P(\nu,x,\pd{x}) \mid P(s,x,\pd{x})\in J_f\}$.
Then $\Osc_X[f^{-1}] \simeq \Dsc_X\otimes_{A_n}(A_n/J_f(\nu))$.  
\item
Compute $\Osc_X[f^{-1}]\otimes_{\Osc_X}\Msc$, which is obtained as 
the output of Algorithm \ref{algorithm:tensor} 
with $A_n/J_f(\nu)$ and $A_n^r/N$ as input and $i=0$.
\end{enumerate}
\end{algorithm}

\section{Algebraic local cohomology groups}
\label{section:localCohomology}
\setcounter{equation}{0}

Let $f_1,\dots,f_d \in K[x]$ be arbitrary polynomials and put
$Y := \{x \in X \mid f_1(x) = \dots = f_d(x) = 0\}$ with $X := K^n$.
Let $\Msc$ be a coherent $\Dsc_X$-module.
Our purpose is to compute the algebraic local cohomology groups
$\Hsc^i_{[Y]}(\Msc)$ with support $Y$ defined by Grothendieck 
as $\Dsc_X$-modules.  
Recall that $\Hsc^i_{[Y]}(\Msc)$ is defined as the $k$-th derived functor 
of the functor 
$$
\Gamma_{[Y]}(\Msc) := \lim_\rightarrow
\Homsc_{\Osc_X}(\Osc_X/\Jsc_Y^m; \Msc),
$$
where $\Jsc_Y$ is the defining ideal of $Y$ and the inductive limit is 
taken as $m$ tends to infinity. 
Note that if $\Msc$ is a holonomic $\Dsc_X$-module, then so is 
$\Hsc^i_{[Y]}(\Msc)$ (cf.\ \cite[Theorem 1.4]{HolonomicII}).  

Put $\widetilde X := K^d \times X$ and identify $X$ with the linear subvariety
$\{0\} \times X$ of $\widetilde X$. We set
$$
Z := \{(t,x) \in \widetilde X \mid t_i = f_i(x)\,\,(i=1,\dots,d)\}.
$$
Then $\Bsc_{[Z]} := \Hsc^d_{[Z]}(\Osc_{\widetilde X})$ is isomorphic to 
$\Dsc_{\widetilde X}/\Isc$, where $\Isc$ is the left ideal generated by 
\begin{equation} \label{eq:BZ}
t_j - f_j(x) \quad (j=1,\dots,d),\qquad
\pd{x_i} + \sum_{j=1}^d\frac{\partial f_j}{\partial x_i}\pd{t_j}
\quad (i=1,\dots,n).
\end{equation}
Let $\pi : \widetilde X \longrightarrow X$ be the projection and put
$$
\Delta':= \{(t,x,y) \in \widetilde X \times X \mid x = y\}.
$$
Then we can identify $\widetilde X$ with $\Delta'$ by $\pi$.   
In the same way as \cite[Lemma 6.3]{OakuAdvance} we get

\begin{lemma}
$$
\Torsc_i^{\pi^{-1}\Osc_X}(\Bsc_{[Z]},\pi^{-1}\Msc) \simeq 
\left\{\begin{array}{ll} 
\Hsc^0((\Bsc_{[Z]}\hat\otimes \Msc)_{\Delta'}^\bullet) & (i=0) \\
0 & (i \neq 0)
\end{array}\right.
$$
with
$$
\Bsc_{[Z]}\hat\otimes\Msc := \Dsc_{\widetilde X \times X}
\otimes_{p_1^{-1}\Dsc_{\widetilde X}\otimes p_2^{-1}\Dsc_X}
(p_1^{-1}\Bsc_{[Z]}\otimes_K p_2^{-1}\Msc),
$$
where $p_1$ and $p_2$ are the projections of $\widetilde X \times X$ to 
$\widetilde X$ and to $X$ respectively.  
\end{lemma}

In fact, this lemma follows from the fact that $\Delta'$ is 
non-characteristic with respect to $\Bsc_{[Z]}\hat\otimes\Msc$.  
The proof of \cite[Theorem 6.1]{OakuAdvance} yields the following: 

\begin{proposition} \label{prop:localCohomology}
For any coherent $\Dsc_X$-module $\Msc$, we have an isomorphism
$$
\Hsc^i_{[Y]}(\Msc) \simeq 
\Hsc^{i-d}((\Bsc_{[Z]}\otimes_{\pi^{-1}\Osc_X}\pi^{-1}\Msc)
_{\{0\}\times X}^\bullet)
$$
as left $\Dsc_X$-module for any $i \geq 0$.
\end{proposition}

\begin{algorithm} \label{algorithm:localCohomology}
\rm
(Algebraic local cohomology groups $\Hsc_{[Y]}^i(\Msc)$)

\noindent
Input: Polynomials $f_1,\dots,f_d \in K[x]$ and 
a holonomic $\Dsc_X$-module $\Msc = \Dsc_X\otimes_{A_n}(A_n^r/N)$
with an $A_n$-submodule $N$ of $A_n^r$ generated by $G$.  
\\
Output: $\Hsc^i_{[Y]}(\Msc) = \Lsc_{i-d}$ for $i=0,\dots,d$ 
with $Y := \{x \in K^n \mid f_1(x) = \dots = f_d(x) = 0\}$.  

\begin{enumerate}
\item
Let $I$ be the left $A_{d+2n}$-submodule of $A_{d+2n}^r$ generated by
\begin{eqnarray*}
\G_1 & := &
\{ (t_j - f_j(x))e_k \mid 1 \leq j \leq d, \,\, 1 \leq k \leq r\}\\
&\cup&
(\pd{x_i} + \sum_{j=1}^d\frac{\partial f_j}{\partial x_i}\pd{t_j})e_k 
\mid 1 \leq i \leq n,\,\, 1 \leq k \leq r\}\\
&\cup&
\{P(y,\pd{y}) \mid P \in G\},
\end{eqnarray*}
where $e_1,\dots,e_r$ are the unit vectors of $A_n^r$.  
\item
Apply substitution $y_i = x_i + z_i$ for $i=1,\dots,n$ to $\G_1$ and 
let the result be $\G_2$.  Let $J_2$ be the submodule of $A_{d+2n}^r$ 
generated by $\G_2$.  
\item 
Compute the $0$-th cohomology of the restriction of $A_{d+2n}^r/J_2$ 
to $\{(t,x,z) \mid z = 0\}$ 
in the form $\Dsc_{\widetilde X}\otimes_{A_{d+n}}(A_{d+n}^l/J_3)$     
by Algorithm \ref{algorithm:restriction}.  Here 
we can assume $k_0 = k_1 = 0$ skipping the steps 1--4 of Algorithm 
\ref{algorithm:restriction}. 
\item
Compute 
$\Lsc_i := 
\Hsc^i((\Dsc_{\widetilde X}\otimes_{A_{d+n}}(A_{d+n}^l/J_3))_
{\{0\}\times X}^\bullet)$ for $-d \leq i \leq 0$.  
\end{enumerate}
\end{algorithm}

Finally, assume that $Y$ is non-singular and let $\iota:Y \rightarrow X$ 
be the embedding.  
Then for a coherent $\Dsc_Y$-module $\Nsc$, 
$$
\iota_+ \Nsc := \Dsc_{X\leftarrow Y}\otimes_{\Dsc_Y}\Nsc
$$
is a coherent $\Dsc_X$-module with support in $Y$; here we put
$\Dsc_{X\leftarrow Y} := \Dsc_X/(\Dsc_Xf_1 + \cdots + \Dsc_Xf_d)$, 
which has a structure of $(\Dsc_X,\Dsc_Y)$-bimodule.  
Moreover, the functor $\iota_+$ gives an equivalence between the 
category of coherent $\Dsc_Y$-modules and that of coherent $\Dsc_X$-modules
supported by $Y$ (\cite[Proposition 4.2]{HolonomicII}).
In terms of this equivalence, we can compute the cohomology groups of 
the restriction of a holonomic $\Dsc_X$-module $\Msc$ to $Y$ 
by using Algorithm \ref{algorithm:localCohomology} and 
the following isomorphism: 

\begin{proposition} {\rm (\cite[Proposition 4.3]{HolonomicII})}
$$
\iota_+ \Hsc^i(\Msc_Y^\bullet) = \Hsc^{i+d}_{[Y]}(\Msc).
$$
\end{proposition}

\begin{example}\rm
Put $K := \C^3 \ni (x,y,z)$ and 
$$
\Msc := \Dsc_X/(\Dsc_X \pd{x} + \Dsc_X \pd{y} + \Dsc_X (z^3\pd{z}+z)).
$$
Then the local cohomology groups of $\Msc$ with support 
$Y := \{(x,y,z) \in X \mid xz = yz = 0\}$ are given by
$$
\Hsc^i_{[Y]}(\Msc) = \left\{ \begin{array}{ll}
\Dsc_X/(\Dsc_X x + \Dsc_X y + \Dsc_X(z^2\pd{z} + 2z + 1)) & (i=2),\\
0 & (i=1),\\
\Dsc_X/(\Dsc_X\pd{x}^2 + \Dsc_X(x\pd{x}-1) + \Dsc_X \pd{y} + \Dsc_X z)
& (i=0).
\end{array}\right.
$$
Moreover, $\Hsc^0_{[Y]}(\Msc) = \Dsc_Xu$ is also isomorphic to 
$$
\Dsc_X/(\Dsc_X \pd{x} + \Dsc_X \pd{y} + \Dsc_X z) = \Dsc_Xv
$$
by the correspondence $u \mapsto \pd{x}v$ and $v \mapsto xu$, where 
$u$ and $v$ denote the residue classes of $1 \in \Dsc_X$ in respective
modules.  

The localization of $\Msc$ by $z$ is given by
$$
\Msc[1/z] = \Dsc_X/(\Dsc_X\pd{x} + \Dsc_X\pd{y} + \Dsc_X(z^2\pd{z}+z+1)).  
$$
Note that the methods of \cite{OakuAdvance} or \cite{Walther} cannot be
appied since $z: \Msc \rightarrow \Msc$ is not injective. 
\end{example}

\section{Functors in the analytic category}
\label{section:analytic}
\setcounter{equation}{0}

The functors studied in the preceding sections have analytic counterparts
(cf.\ \cite{HolonomicII},\cite{MebkhoutBook}).
Throughout this section, we assume $K= \C$ and denote by $\an{O}{X}$ and 
$\an{D}{X}$ the sheaves on $X$ of holomorphic functions and of rings of 
differential operators with holomorphic coefficients respectively.  
For a left $\Dsc_X$-module $\Msc = \Dsc_X^r/\Nsc$, we put 
$\an{M}{} := \Dsc_X^{\rm an}\otimes_{\Dsc_X}\Msc$.  
For $k \in \Z$, put
$$
F_Y^k(\an{D}{X}) := \{ P \in \an{D}{X}|_Y \mid P (\an{J}{Y})^j 
\in (\an{J}{Y})^{j-k} \quad\mbox{for any $j\geq k$} \},
$$
where $\an{J}{Y}$ is the defining ideal of $Y$ in $\an{O}{X}$.   
Then for $\mvec = (m_1,\dots,m_r) \in \Z^r$, 
the $F_Y[\mvec]$-filtrations are defined by
\begin{eqnarray*}
F_Y^k[\mvec]((\an{D}{X})^r) &:=& \bigoplus_{i=1}^r F_Y^{k-m_i}(\an{D}{X}),\\
F_Y^k[\mvec](\an{M}{}) &:=& F_Y^{k-m_1}(\an{D}{X})u_1 + \cdots + 
                        F_Y^{k-m_r}(\an{D}{X})u_r,
\end{eqnarray*} 
where $u_1,\dots,u_r$ are the residue classes of the unit vectors of 
$(\an{D}{X})^r$.  The graded modules are defined in the same way as in 
Section \ref{section:V-filtration}.
Put $\an{D}{Y\rightarrow X} := \an{O}{Y}\otimes_{\an{O}{X}}\an{D}{X}$, 
which is a $(\an{D}{Y},\an{D}{X})$-bimodule. 
Then the restriction of $\an{M}{}$ to $Y$ is defined by
$$
(\an{M}{})_Y^\bullet := 
\an{D}{Y\rightarrow X}\otimes^L_{\an{D}{X}}\an{M}{}. 
$$
The $b$-function of $\an{M}{}$ along $Y$ at $p\in Y$ is defined to be the
generator of the ideal 
$$
\{ b(\theta) \in \C[\theta] \mid 
b(\vartheta)\gr_Y^0[\mvec](\an{M}{}) = 0 \},
$$
where $\vartheta$ is defined as in Section \ref{section:b-function}.

\begin{lemma} \label{lemma:comparison}
Let $\Nsc$ be a coherent left $\Dsc_X$-submodule of $\Dsc_X^r$.  
Define the $F_Y[\mvec]$-filtrations on $\Nsc$ and on $\an{N}{}$ as in 
Section \ref{section:V-filtration} with a shift vector $\mvec$.
Then we have
\begin{eqnarray*}
\gr_Y^0[\mvec](\an{N}{}) 
&=& \gr_Y^0(\an{D}{X})\otimes_{\gr_Y^0(\Dsc_X)}\gr_Y^0[\mvec](\Nsc) \\
&=& \an{D}{Y}\otimes_{\Dsc_Y}\gr_Y^0[\mvec](\Nsc).
\end{eqnarray*}
\end{lemma}

\begin{proof}
We can prove the first equality by the same method (considering 
syzygies in the graded module) as \cite[Theorem 3.16]{OakuBfunction}
(cf.\ also \cite[Lemma 1.1.2]{ACG}), where the case of $r=d=1$ 
is treated; the argument applies to this case with trivial modifications.
The second equality follows from 
$$
\gr_Y^0(\Dsc_X) = \Dsc_Y[t_1\pd{t_1},\dots,t_d\pd{t_d}],\qquad
\gr_Y^0(\an{D}{X}) = \an{D}{Y}[t_1\pd{t_1},\dots,t_d\pd{t_d}].
$$
\end{proof}

\begin{proposition}
The $b$-function $b^{\rm an}(\theta,p)$ of $\an{M}{}$ and the 
$b$-function $b(\theta,p)$ of $\Msc$ with the same shift vector 
$\mvec$ coincide for any $p \in Y$. 
\end{proposition}

\begin{proof}
This follows from Lemma \ref{lemma:comparison} and the faithful flatness of
$\an{O}{Y}$ over $\Osc_Y$ (cf.\ \cite[Lemma 4.4]{OakuBfunction}).
\end{proof}

\begin{proposition} \label{prop:comparison-restriction}
We have for any $i\in \Z$, 
$$
\Hsc^i((\an{M}{})_Y^\bullet) = 
\an{D}{Y}\otimes_{\Dsc_Y}\Hsc^i(\Msc_Y^\bullet).
$$
\end{proposition}

\begin{proof}
Since we can regard
$$
\an{D}{Y\rightarrow X} = 
\{ \sum_{\nu,\beta}a_{\nu\beta}(x)\pd{t}^\nu\pd{x}^\beta \in \an{D}{X} \mid 
a_{\nu\beta}(x) \in \an{O}{Y}\},
$$
we have an isomorphism 
$\an{D}{Y\rightarrow X} \simeq \an{D}{Y}\otimes_{\Dsc_Y}\Dsc_{Y\rightarrow X}$ 
as $(\an{D}{Y},\Dsc_X)$-bimodules. 
Combining this with the flatness of $\an{D}{X}$ over $\Dsc_X$, 
and that of $\Dsc_{Y\rightarrow X}$ over $\Dsc_Y$, we get
\begin{eqnarray*}
(\an{M}{})_Y^\bullet &=&  \an{D}{Y\rightarrow X} \otimes^L_{\an{D}{X}}
(\an{D}{X}\otimes_{\Dsc_X}\Msc) \\
&=& \an{D}{Y\rightarrow X} \otimes^L_{\Dsc_X}\Msc \\
&=& (\an{D}{Y} \otimes_{\Dsc_Y} \Dsc_{Y\rightarrow X})\otimes^L_{\Dsc_X}\Msc \\
&=& \an{D}{Y}\otimes_{\Dsc_Y}(\Dsc_{Y\rightarrow X}\otimes^L_{\Dsc_X}\Msc) \\
&=& \an{D}{Y}\otimes_{\Dsc_Y}\Msc_Y^\bullet. \\
\end{eqnarray*}
This implies the assertion since $\an{D}{Y}$ is flat over $\Dsc_Y$. 
\end{proof}

\begin{proposition} \label{prop:comparison-tensor}
For $\Dsc_X$-modules $\Msc$ and $\Nsc$, we have 
$$
\Torsc^{\an{O}{X}}_i(\an{M}{},\an{N}{}) =
 \an{D}{X}\otimes_{\Dsc_X}\Torsc^{\Osc_X}_i(\Msc,\Nsc).
$$
\end{proposition}

\begin{proof}
The assertion follows from Proposition \ref{prop:comparison-restriction}
and \cite[Proposition 4.7]{HolonomicII}.  
\end{proof}

For an algebraic set $Y$ of $X$, the algebraic local cohomology groups 
of $\an{M}{}$ are defined to be the derived functors of the functor
$$
\Gamma_{[Y]}(\an{M}{}) := \lim_{m\rightarrow\infty}
\Homsc_{\an{O}{X}}(\an{O}{X}/(\an{J}{Y})^m; \an{M}{}).
$$

\begin{proposition}
For a left $\Dsc_X$-module $\Msc$, a polynomial $f\in \C[x]$, and 
an algebraic set $Y$, we have
\begin{eqnarray*}
\an{M}{}[f^{-1}] &=& \an{D}{X}\otimes_{\Dsc_X}\Msc[f^{-1}],\\
\Hsc_{[Y]}^i(\an{M}{}) &=& 
\an{D}{X}\otimes_{\Dsc_X}\Hsc_{[Y]}^i(\Msc).
\end{eqnarray*}
\end{proposition}

\begin{proof}
The first equality is an immediate consequence of Proposition 
\ref{prop:comparison-tensor}. 
The second equality follows from Proposition \ref{prop:comparison-restriction},
and Proposition \ref{prop:localCohomology} together with 
its analytic counterpart.  
\end{proof}

\section{Homogenized Weyl algebra and Schreyer's method for free resolution}
\label{section:Schreyer}
\setcounter{equation}{0}

In this section, we work in a framework more general than is needed in the
preceding sections.  
Let $A_n = K[x]\langle\partial\rangle$ be the Weyl algebra over a field 
$K$ of characteristic zero with $x = (x_1,\dots,x_n)$,
$\partial = (\pd{1},\dots,\pd{n})$, $\pd{i} = \partial/\partial x_i$.   
We introduce a weight vector 
$w = (w_1,\dots,w_n;w_{n+1},\dots,w_{2n})\in \Z^{2n}\setminus\{0\}$ 
that satisfies  $w_i + w_{n+i} \geq 0$ for $i=1,\dots,n$.  
For each integer $\nu\in\Z$, we put
$$
F_w^\nu(A_n) := \{P = \sum_{\alpha,\beta\in\N^n}
a_{\alpha\beta}x^\alpha\partial^\beta \in A_n \mid
a_{\alpha\beta}= 0\, 
\mbox{if } \sum_{i=1}^n w_i\alpha_i + \sum_{i=1}^n w_{n+i}\beta_i > \nu \}, 
$$
where $a_{\alpha\beta}\in K$ and the sum with respect to 
$\alpha,\beta \in \N^n$ is finite.  
For a nonzero $P \in A_n$, let $\ord_w(P)$ be the minimum integer $k$ 
such that $P \in F_w^k(A_n)$.  
It is easy to see that $\ord_w(PQ) = \ord_w(P)+\ord_w(Q)$ holds for 
nonzero $P,Q \in A_n$.
More generally, for a shift vector $\mvec = (m_1,\dots,m_r) \in \Z^r$, 
we define a filtration 
$F_w[\mvec]$ of $A_n^r$ by
$$
F_w^k[\mvec](A_n^r) := \bigoplus_{i=1}^r  F_w^{k-m_i}(A_n)e_i,
$$
where $e_1,\dots,e_r$ are the canonical generators of $A_n^r$.  
For a nonzero $P \in A_n^r$, we put
$\ord_w[\mvec](P) := \min\{k\in\Z \mid P \in F_w^k[\mvec](A_n^r)\}$. 

Now we introduce the homogenized Weyl algebra, which was introduced from
the second version (1994) of {\tt kan/sm1} \cite{Kan-sm1}: 

\begin{definition}
\rm
Let $A_n^{(h)}$ be the algebra over $K$ generated by 
$h$, $x = (x_1,\dots,x_n)$, and $\partial = (\pd{1},\dots,\pd{n})$ which
satisfy the relations
$$ \matrix{
& x_ix_j - x_jx_i = 0, & \partial_i\partial_j - \partial_j\partial_i = 0,  
& x_i\partial_j - \partial_jx_i = -\delta_{ij}h^2, & \cr
& hx_i - x_ih = 0, & h\partial_i-\partial_i h = 0 & (i,j = 1,\dots,n), \cr}
$$
where $\delta_{ii} = 1$ and $\delta_{ij}=0$ if $i \neq j$.
We call $A_n^{(h)}$ the {\em homogenized Weyl algebra}.  
The substitution $h=1$ defines a $K$-algebra homomorphism 
$$ \rho \,:\,\, A_n^{(h)} \ni P \longmapsto P|_{h=1} \in A_n. $$
\end{definition}

An element $P$ of $A^{(h)}$ is uniquely expressed as a finite sum
$$
P = \sum_{\lambda\in\N,\alpha,\beta\in\N^n}a_{\lambda\alpha\beta}
h^\lambda x^\alpha \partial^\beta
$$
with $a_{\lambda\alpha\beta}\in K$. 
The total degree of $P \in A^{(h)}$ is defined by
$$
\deg(P) := \max\{ \lambda + |\alpha| + |\beta| \mid a_{\lambda\alpha\beta}
\neq 0 \}
$$
if $P \neq 0$ and $\deg(P) = -\infty$ if $P=0$.  

Now let us take another vector $\nvec = (n_1,\dots,n_r)\in\Z^r$, 
which describes the shift with respect to the total degree.  
For $P = (P_1,\dots,P_r) \in (A^{(h)}_n)^r$, we put
$$ 
\deg[\nvec](P) := \max_{1\leq i \leq r} (\deg(P_i) + n_i). 
$$

\begin{definition}
\rm
\begin{enumerate}
\item
An element $P = \sum_{i=1}^r\sum_{(\lambda,\alpha,\beta)\in L} 
a_{\lambda\alpha\beta i}h^\lambda x^\alpha\partial^\beta e_i$ of 
$(A_n^{(h)})^r$ 
is said to be $h[\nvec]$-homogeneous if there exists $k \in \Z$ so that
$a_{\lambda\alpha\beta i}=0$ unless $\lambda + |\alpha|+|\beta|+n_i = k$.
\item
For $P = \sum_{i=1}^r\sum_{(\alpha,\beta)\in L} 
a_{\alpha\beta i}x^\alpha\partial^\beta e_i \in A_n^r$, we define 
the $h[\nvec]$-homogenization 
$h[\nvec](P)\in (A^{(h)}_n)^r$ by
$$
h[\nvec](P) := \sum_{i=1}^r\sum_{(\alpha,\beta)\in L} 
a_{\alpha\beta i}h^{k-|\alpha|-|\beta|-n_i}x^\alpha\partial^\beta e_i
$$
with $k := \max\{|\alpha|+|\beta|+n_i \mid a_{\alpha\beta i} \neq 0\}$.
This is $h[\nvec]$-homogeneous.
\end{enumerate}
\end{definition}

If $\nvec$ is the zero vector, we denote $h[\nvec](P)$ simply by $h(P)$.

\begin{lemma}
For $P\in A_n^r$ and $Q \in A_n$,  we have 
$\rho(h[\nvec](P)) = P$ and $h[\nvec](QP) = h(Q)h[\nvec](P)$.
\end{lemma}

\begin{definition} \label{def:adapted-order}
\rm 
Let $\prec$ be a monomial order 
(i.e.\ an order satisfying (\ref{eq:monomial-order})) 
on $L \times \Set{r}$ with $L := \N^{2n}$.  
We denote by $\lexp_\prec(P)$ the leading exponent of $P \in A_n^r$ 
with respect to $\prec$.  
Then $\prec$ is said to be {\em adapted to} the filtration $F_w[\mvec]$ if 
$\langle w,\alpha \rangle + m_i < \langle w,\beta\rangle + m_j$ implies 
$(\alpha,i) \prec (\beta,j)$ for $\alpha,\beta \in L$ and $i,j\in\Set{r}$,
and if $\lexp_\prec(e_i) \prec \lexp_\prec(x_j\pd{j}e_i)$ 
for any $1\leq i \leq r$ and $1 \leq j \leq n$; 
here we write $\langle w,\alpha\rangle = \sum_{i=1}^{2n}w_i\alpha_i$ for
$\alpha = (\alpha_1,\dots,\alpha_{2n})$.  
\end{definition}

We fix a monomial order $\prec$ on $L\times\Set{r}$ that is adapted to the 
$F_w[\mvec]$-filtration.  
Then we define an order $\prec_h$ on 
$\N\times L\times\Set{r}$ by
\begin{eqnarray*}
& & 
\mbox{
$(\lambda,\alpha,i) \prec_{h[\nvec]} 
(\mu,\beta,j)$ if and only if}\\
& &
\mbox{
$\lambda + |\alpha|+n_i < \mu + |\beta|+n_j$ 
or else}\\
& &
\mbox{
$\lambda + |\alpha|+n_i = \mu + |\beta|+n_j, \,\,
(\lambda,\alpha,i) \prec (\mu,\beta,j)$,}
\end{eqnarray*}
for $\lambda,\lambda' \in \N, \alpha,\beta \in L$ and 
$i,j \in \{1,...,r\}$.  
Then it is easy to see that $\prec_{h[\nvec]}$ is a well-order. 
For a nonzero element 
\begin{equation} \label{eq:PinAhr}
P = \sum_{i=1}^r \sum_{\lambda,\alpha,\beta}a_{\lambda\alpha\beta i}
h^\lambda x^\alpha \partial^\beta e_i
\end{equation}
of $(A^{(h)}_n)^r$, its leading exponent 
$(\lambda_0,\alpha_0,\beta_0, i_0) = \lexp_{h[\nvec]}(P) 
\in \N\times L \times \Set{r}$ is defined as the maximum element of 
$\{(\lambda,\alpha,\beta,i) \mid a_{\lambda\alpha\beta i} \neq 0 \}$
with respect to $\prec_{h[\nvec]}$.  
Then the leading position $\lp_{h[\nvec]}(P)$ and the leading coefficient 
$\lcoef_{h[\nvec]}(P)$ are defined to be $i_0$ and 
$a_{\lambda_0,\alpha_0,\beta_0,i_0}$ respectively.
We denote them simply by $\lexp(P)$, $\lp(P)$, and $\lcoef(P)$ if 
there is no fear of confusion.  
The following lemmas follow easily from Definition \ref{def:adapted-order}  
and the definitions of $A^{(h)}_n$ and $\prec_{h[\nvec]}$:

\begin{lemma} \label{lemma:lexp}
For $P \in (A^{(h)}_n)^r$ and $Q \in A^{(h)}_n$, we have
$\lexp(QP) = \lexp(Qe_k)+\lexp(P)$ 
with $k = \lp(P)$.   
\end{lemma}

\begin{lemma}
If $P \in (A^{(h)}_n)^r$ is $h[\nvec]$-homogeneous 
and $Q \in A^{(h)}_n$ is $h[0]$-homogeneous, then
$QP$ is $h[\nvec]$-homogeneous. 
\end{lemma}

\begin{lemma}
Let $\varpi : \N \times L \times \Set{r} \rightarrow L \times \Set{r}$ 
be the projection.  
Then $\lexp_\prec(\rho(P)) = \varpi(\lexp_{h[\nvec]}(P))$ holds 
if $P$ is $h[\nvec]$-homogeneous. 
\end{lemma}

In view of the above lemmas, we can define the notion of Gr\"obner basis 
in the homogenized Weyl algebra and can employ the Buchberger algorithm, 
which preserves the $h[\nvec]$-homogeneity: 

\begin{definition}
\rm
Let $N$ be a left $A^{(h)}_n$-submodule of $(A^{(h)}_n)^r$.  
Then a finite subset $G$ of $N$ is called a Gr\"obner basis of $N$ 
with respect to $\prec_{h[\nvec]}$ if
$$ 
E(N):= \{\lexp(P) \mid P \in N\setminus\{0\}\} 
= \bigcup_{P \in \G} (\lexp(P) + L).
$$
\end{definition}
Note that $G$ generates $N$ if $G$ is a Gr\"obner basis of $N$ since 
$\prec_{h[\nvec]}$ is a well-order. 

\begin{proposition} \label{prop:h-division}
Let $N$ be a left $A_n$-submodule of $A_n^r$ generated by $P_1,\dots,P_k$.  
Let $h[\nvec](N)$ be the left $A^{(h)}_n$-submodule of $(A^{(h)}_n)^r$ 
generated by $h[\nvec](P_i)$ for $i=1,\dots,k$.  
Let $G = \{Q_1,\dots,Q_s\}$ be a Gr\"obner basis of $h[\nvec](N)$ 
with respect to $\prec_{h[\nvec]}$.  
Then for any $P \in N$, there exist $U_j \in A_n$ such that 
$P = U_1\rho(Q_1) + \cdots + U_s\rho(Q_s)$ and 
$\ord_w[\mvec](U_j\rho(Q_j)) \leq \ord_w[\mvec](P)$ for $j=1,...,s$.
\end{proposition}

\begin{proof}
There exists $\nu\in\N$ such that $h^\nu h[\nvec](P)$ belongs to 
$h[\nvec](N)$.  
By the division algorithm in $(A^{(h)}_n)^r$, we can find 
$h[\nvec]$-homogeneous $U_1,\dots,U_s \in (A^{(h)})^r$ such that 
$h^\nu h[\nvec](P) = U_1Q_1 + \cdots + U_sQ_s$  and
$\lexp(U_kQ_k) \preceq_{h[\nvec]} \lexp(h^\nu h[\nvec](P))$ for $k=1,\dots,s$. 
Applying the ring homomorphism $\rho$, we get
$P = \rho(U_1)\rho(Q_1) + \cdots + \rho(U_s)\rho(Q_s)$ 
and $\lexp_\prec(\rho(U_k)\rho(Q_k)) = \varpi(\lexp(U_kQ_k))
\preceq \varpi(\lexp(h^\nu h[\nvec](P))) = \lexp_\prec(P)$.
Since $U_i$ and $Q_i$ are $h[\nvec]$-homogeneous, this implies
$\ord_w[\mvec](\rho(U_iQ_i)) \leq \ord_w[\mvec](P)$.  
\end{proof}

We use the same notation as in the preceding proposition. 
Put $\Lambda := \{(i,j) \mid 1 \leq i < j \leq s,\,\,\lp(P_i) = \lp(P_j)\}$.  
Then for $(i,j) \in \Lambda$, let $S_{ij}, S_{ji} \in A^{(h)}_n$ 
be monomials such that 
$$
\lexp(S_{ji}P_i) = \lexp(S_{ij}P_j) = \lexp(P_i) \vee \lexp(P_j),
\quad
\lcoef(S_{ji}P_i) = \lcoef(S_{ij}P_j).
$$
By the Buchberger algorithm, there exist $h[\nvec]$-homogeneous 
$U_{ijk} \in A^{(h)}_n$ so that we have
$$
S_{ji}P_i - S_{ij}P_j = \sum_{k=1}^s U_{ijk}P_k
$$
and either $U_{ijk} \neq 0$ or 
$$
\lexp(U_{ijk}P_k) \prec_{h[\nvec]} \lexp(P_i) \vee \lexp(P_j)
$$ 
for each $k = 1,...,s$.  
The following is an analogue of F.O.~Schreyer's theorem for the syzygies in 
the polynomial ring: 

\begin{theorem} \label{th:Schreyer}
In the notation above, 
let $\prec'$ be the order on $\N\times L \times \Set{s}$ defined by
\begin{eqnarray*}
& & (\alpha,\mu) \prec' (\beta,\nu) \quad\mbox{if and only if}\quad
\lexp(P_\mu)+\alpha \prec_{h[\nvec]} \lexp(P_\nu)+\beta \\
& & \mbox{or else}\quad
\lexp(P_\mu)+\alpha = \lexp(P_\nu)+\beta 
\mbox{ and } \mu > \nu
\end{eqnarray*}
for $\alpha,\beta\in \N\times L$ and $\mu,\nu\in\Set{s}$.  
Put 
$$
\mvec' := (\ord_w[\mvec](P_1),\dots,\ord_w[\mvec](P_s)),\quad
\nvec'= (\deg[\nvec](P_1),\dots,\deg[\nvec](P_s)).
$$ 
Then $\prec'$ is a well-order adapted to the $F_w[\mvec']$-filtration 
and called the Schreyer order induced by $\prec_{h[\nvec]}$.
\begin{enumerate}
\item For $(i,j)\in\Lambda$, 
$$
V_{ij} := (0,\dots,\stackrel{(i)}{S_{ji}},\dots,
\stackrel{(j)}{-S_{ij}},\dots,0) - (U_{ij1},\dots,U_{ijs})
$$
is $h[\nvec']$-homogeneous and 
$\{V_{ij} \mid (i,j) \in \Lambda \}$ is a Gr\"obner basis 
with respect to $\prec'$ of the module
$$
{\rm Syz}(P_1,\dots,P_s) 
:= \{(U_1,...,U_s) \in (A^{(h)}_n)^s \mid U_1P_1 + \cdots + U_sP_s = 0\}. 
$$
\item
Put  
$$
{\rm Syz}(\rho(P_1),\dots,\rho(P_s)) 
:= \{(U_1,...,U_s) \in A_n^s \mid U_1\rho(P_1) + \cdots + U_s\rho(P_s) = 0\}. 
$$
Then for any $P \in {\rm Syz}(\rho(P_1),\dots,\rho(P_s)) \cap 
F_w^\nu[\mvec'](A_n^s)$ with $\nu \in \Z$, there exist
$Q_{ij}\in A_n$ such that
$P = \sum_{(i,j)\in\Lambda}Q_{ij}\rho(V_{ij})$ and that
$ Q_{ij}\rho(V_{ij}) \in F_w^\nu[\mvec'](A_n^s)$.
\end{enumerate}
\end{theorem}

\begin{proof}
For a nonzero $P \in (A^{(h)}_n)^r$ of the form (\ref{eq:PinAhr}), we 
define its initial term by
$$
\init(P) := a_{\lambda\alpha\beta i}h^\lambda x^\alpha
\xi^\beta e_i \,\,\in K[h,x,\xi]^r
$$
with $(\lambda,\alpha,\beta,i) = \lexp(P)$,  
where $\xi = (\xi_1,\dots,\xi_n)$ are commutative indeterminates. 
Let $s_{ij}$ be the monomial in $K[h,x,\xi]$ 
obtained by substituting $\xi$ for $\partial$ in $S_{ij}$.  
Then we have $s_{ji}\init(P_i) - s_{ij}\init(P_j) = 0$
for $(i,j) \in \Lambda$.  
Now suppose $U = (U_1,\dots,U_s) \in {\rm Syz}(P_1,\dots,P_s)$ and put 
$(\lambda_0,\alpha_0,\beta_0,i_0) := \max_{1\leq j \leq s}\lexp(U_jP_j)$.  
Define $u_j \in K[h,x,\xi]$ by 
$$
u_je_{i_0} = \init(U_je_{i_0}) \quad\mbox{if 
$\lexp(U_jP_j) = (\lambda_0,\alpha_0,\beta_0,i_0)$},
$$
and put $u_j = 0$ otherwise.
Then $\sum_{j=1}^s u_j\init(P_j) = 0$ holds. 
By the definition of $\prec'$ we have
\begin{equation} \label{eq:Vij}
\lexp_{\prec'}(V_{ij}) = (\lexp(s_{ji}),i),\quad
\lexp_{\prec'}(U) = \lexp_{\prec'}((u_1,\dots,u_s)).
\end{equation}  
On the other hand,
since $(0,\dots,s_{ji},0,\dots,-s_{ij},0,\dots,0)$ for $(i,j)\in\Lambda$ 
constitute a Gr\"obner basis with respect to $\prec'$ of the syzygies 
on $\init(P_1),\dots,\init(P_s)$ by virtue of Schreyer's theorem for 
the polynomial ring (cf.\ \cite[Theorem 15.10]{EisenbudBook}), we know that
$$
\lexp_{\prec'}((u_1,\dots,u_s))\in 
\bigcup_{(i,j)\in\Lambda} ((\lexp(s_{ji}),i) + (\N \times L)).
$$ 
This completes the proof of the first assertion in view of (\ref{eq:Vij}).  

The second assertion follows from the first and Proposition 
\ref{prop:h-division}
\end{proof}

Let $N$ be a left $A_n$-submodule of $A_n^r$ generated by $P_1,\dots,P_k$. 
Let $h(N)$ be the left $A^{(h)}_n$-submodule of $(A^{(h)}_n)^r$ 
generated by $h(P_1),\dots,h(P_k)$ (the homogenizations with $\nvec = 0$).  
Starting with $h(N)$ and $\nvec = \mvec = 0$, 
apply the first part of Theorem \ref{th:Schreyer} repeatedely.  
Then we get an exact sequence 
\begin{equation} \label{eq:Ah-resolution}
\cdots 
\stackrel{\psi_3}{\longrightarrow} (A^{(h)}_n)^{r_2}
\stackrel{\psi_2}{\longrightarrow} (A^{(h)}_n)^{r_1}
\stackrel{\psi_1}{\longrightarrow} (A^{(h)}_n)^{r_0}
\stackrel{\varphi}{\longrightarrow}(A^{(h)}_n)^{r}/h(N)  \longrightarrow 0
\end{equation}
with $r_0 := r$.  Put  $\mvec_0 = \nvec_0 = 0$ and 
\begin{eqnarray*}
\mvec_i &:=& (\ord_w[\mvec_{i-1}](\rho(\psi_i(1,0,\dots,0))),\dots,
\ord_w[\mvec_{i-1}](\rho(\psi_i(0,\dots,0,1)))), \\
\nvec_i &:=& (\deg[\nvec_{i-1}](\psi_i(1,0,\dots,0)),\dots,
\deg[\nvec_{i-1}](\psi_i(0,\dots,0,1))).
\end{eqnarray*}
 
Applying the homomorphism $\rho$ to (\ref{eq:Ah-resolution}), 
we get an exact sequence
\begin{equation} \label{eq:Fw-resolution}
\cdots 
\stackrel{\rho(\psi_3)}{\longrightarrow} A_n^{r_2}
\stackrel{\rho(\psi_2)}{\longrightarrow} A_n^{r_1}
\stackrel{\rho(\psi_1)}{\longrightarrow} A_n^{r_0}
\stackrel{\rho(\varphi)}{\longrightarrow} A_n^{r}/N  \longrightarrow 0. 
\end{equation}
Moreover, in view of the second part of Theorem \ref{th:Schreyer}, 
the sequence
$$
\cdots 
\stackrel{\rho(\psi_2)}{\longrightarrow} F_w^k[\mvec_1](A_n^{r_1})
\stackrel{\rho(\psi_1)}{\longrightarrow} F_w^k[\mvec_0](A_n^{r_0})
\stackrel{\rho(\varphi)}{\longrightarrow} 
F_w^k[0](A_n^r)/(N \cap F_w^k[0](A_n^r))  \longrightarrow 0
$$
is exact for any $k\in\Z$.  That is, the resolution (\ref{eq:Fw-resolution})
is adapted to the $F_w$-filtration.  
Furthermore, we can prove the following in the same way as its counterpart
in the polynomial ring (\cite[Corollary 15.11]{EisenbudBook})

\begin{theorem} \label{th:finite-resolution}
By arranging the Gr\"obner bases appropriately, we can construct the 
free resolution (\ref{eq:Ah-resolution}) so that $\psi_{2n+2}=0$.
\end{theorem}

It is known that there exists a shorter projective resolution.
It is an open question to obtain shorter resolutions.

\begin{remark}
\rm 
\begin{enumerate}
\item
If each element of the weight $w$ is non-negative and the order $\prec$ 
adapted to $F_w[\mvec]$ is a well-order, 
then the above construction can be done 
directly in $A_n^r$ without the homogenized Weyl algebra.  
\item
If the weight $w$ satisfies $w_i+w_{n+i}=0$ for $i=1,\dots,n$, 
then we can work in $A_n[x_0]^r$ as in Section \ref{section:Gbase} 
instead of the homogenized Weyl algebra.
Then Theorem \ref{th:Schreyer} holds with the $h[\nvec]$-homogenization
replaced by the $F_w[\mvec]$-homogenization defined by 
$$
h[\mvec](P) := \sum_{i=1}^r\sum_{\alpha,\beta\in\N^n} a_{\alpha\beta i}
{x_0}^{\langle w,(\alpha,\beta)\rangle +m_i-k}
x^\alpha\pd{}^\beta e_i
$$
with $k := \min \{\langle w, (\alpha,\beta)\rangle + m_i \mid 
a_{\alpha\beta i} \neq 0\}$ for
$P = \sum_{i=1}^r\sum_{\alpha,\beta\in\N^n} a_{\alpha\beta i}
x^\alpha\pd{}^\beta e_i$.

Theorem \ref{th:finite-resolution} also holds in this case.  
\end{enumerate}
\end{remark}

In the computation of the free resolution described in this section, 
we can employ the method of La Scala and Stillman \cite{LaScala},
which computes the `Schreyer frame' (initial terms of the resolution)
first, then computes the resolution by a selection strategy or in parallel.  
We have implemented this algorithm in {\tt kan/sm1}.  
Most of the examples of this paper have been computed by using 
this implementation.  

\begin{example} \rm
We give an example of the Shreyer resolution and explain the limits
of our method (cf. Example 7.1 in \cite{Walther}) caused by the complexity. 
Let $I$ be the left ideal in
$$ A_9 := 
K \langle x_1, x_2, x_3, x_4, x_5, x_6, 
          \pd{x_1}, \pd{x_2}, \pd{x_3}, \pd{x_4}, \pd{x_5}, \pd{x_6},
             t_1, t_2, t_3, \pd{t_1}, \pd{t_2}, \pd{t_3} \rangle
$$
generated by
\begin{eqnarray*}
 &&x_4 x_2 -x_5 x_1 +t_1  
 , x_4 x_3 -x_6 x_1 +t_2 
 , x_5 x_3 -x_6 x_2 +t_3 \\
 &&  x_5  \pd{t_1}+x_6   \pd{t_2}+\pd{x_1} 
 , -x_4  \pd{t_1}+x_6  \pd{t_3}+\pd{x_2} 
 , -x_4  \pd{t_2}-x_5  \pd{t_3}+\pd{x_3} \\
 && -x_2  \pd{t_1}-x_3  \pd{t_2}+\pd{x_4} 
 , x_1  \pd{t_1}-x_3  \pd{t_3}+\pd{x_5} 
 , x_1  \pd{t_2}+x_2  \pd{t_3}+\pd{x_6} .
\end{eqnarray*}

We want to obtain the restriction of $A_9/I$ along
$t_1 = t_2 = t_3 =0$, but we have not yet succeeded in  
getting the restrictions of all degrees except zero 
by using our algorithm and our implementation; 
we could get only a huge resolution.  
The hardest part in the computation of the restriction is to
compute the quotient of the kernel over the image of the truncated
complex (the step 7 of Algorithm \ref{algorithm:restriction}).
We cannot compute the higher order quotients by our algorithm
and implementations.

We homogenize the ideal with the variable $s$;
\begin{eqnarray*}
 &&x_4 x_2 s-x_5 x_1 s+t_1 
 , x_4 x_3 s-x_6 x_1 s+t_2
 , x_5 x_3 s-x_6 x_2 s+t_3 \\
 &&  x_5 s \pd{t_1}+x_6  s \pd{t_2}+\pd{x_1}
 ,  -x_4 s \pd{t_1}+x_6 s \pd{t_3}+\pd{x_2} 
 ,  -x_4 s \pd{t_2}-x_5 s \pd{t_3}+\pd{x_3} \\
 && -x_2 s \pd{t_1}-x_3 s \pd{t_2}+\pd{x_4} 
 ,  x_1 s \pd{t_1}-x_3 s \pd{t_3}+\pd{x_5} 
 ,  x_1 s \pd{t_2}+x_2 s \pd{t_3}+\pd{x_6} ,
\end{eqnarray*}
which we denote by $I'$.
We compute the Schreyer resolution with the weight matrix
$$\bordermatrix{ &s & x_i & \pd{x_i} & t_i \cr
                 &1 & 0   &  0       &  0  \cr
                 &0 & 1   &  1       &  1  \cr
}. $$
The reduced (non-$h$-homogenized) Gr\"obner basis of $I'$
consists of $55$ elements.
The Betti numbers (i.e., $r_1,r_2,\dots$) of the Schreyer resolution are 
$630, 3329, \ldots$
where the reduced Gr\"obner basis of syzygies by the Schreyer order 
among the $55$-elements
of the reduced Gr\"onber basis of $I'$ consists of $630$ elements.
The syzygy among $3329$ elements by the Schreyer order could not be obtained,
because of the memory exhaustion,
on three MMX Pentium PCs (166 MHz) with 64 Mega bytes of memory 
under the Linux operating  system (version 2.0.30)
with {\tt kan/sm1} (version 2.980129)  and the open sm1 module
for communication between distributed processors \cite{noro-takayama}.
The first Gr\"obner basis of $I'$ can be obtained in two seconds.

The number of Gr\"obner basis is smaller in this case
if we compute in the homogenized Weyl algebra.
We homogenize the ideal $I$ in the homogenized Weyl algebra with the
variable $h$;
\begin{eqnarray*}
  && x_4 x_2-x_5 x_1+t_1 h ,\ 
     x_4 x_3-x_6 x_1+t_2 h , \ 
     x_5 x_3-x_6 x_2+t_3 h ,\\
  && x_5 \pd{t_1}+x_6 \pd{t_2}+h \pd{x_1} ,\ 
     -x_4 \pd{t_1}+x_6 \pd{t_3}+h \pd{x_2} ,\  
     -x_4 \pd{t_2}-x_5 \pd{t_3}+h \pd{x_3} , \\
  && -x_2 \pd{t_1}-x_3 \pd{t_2}+h \pd{x_4} ,\ 
      x_1 \pd{t_1}-x_3 \pd{t_3}+h \pd{x_5} , \ 
      x_1 \pd{t_2}+x_2 \pd{t_3}+h \pd{x_6},
\end{eqnarray*}
which we denote by $I^h$.
We compute the Schreyer resolution with the weight matrix
$$\bordermatrix{  & x_i & \pd{x_i} & t_i & \pd{t_i} \cr
                  & 0   &  0       &  1  &  -1 \cr
                  & 1   &  1       &  1  &  0 \cr
}. $$
The reduced Gr\"obner basis of $I^h$ consists of $44$ elements which 
are less than those of V-homogenized Gr\"obner basis.
The basis can be obtained in 1.5 seconds by the same system.
The Betti numbers of the Schreyer resolution are
$506$, $2422$, $\ldots$.
\end{example}

We conjecture that in the homogenized Weyl algebra the minimal resolution
exists and can be constructed by the algorithm in \cite{LaScala}
since it is a graded algebra, 
but the boundary maps are not adapted to the filtration in general.
For example, 
let us consider the following free resolution
in the homogenized Weyl algebra
$$  0 \longrightarrow A_2^{(h)} \stackrel{\psi_2}{\longrightarrow}
                      (A_2^{(h)})^2 \stackrel{\psi_1}{\longrightarrow}
                      A_2^{(h)} \stackrel{\varphi}{\longrightarrow}
                      A_2^{(h)}/(p_1,p_2) \longrightarrow 0
$$
where
$$ A_2^{(h)} = {\bf C} \langle h, x, y , \partial_x, \partial_y \rangle,
$$
$$ \psi_1( (1,0) ) = p_1, \  \psi_2( (0,1) ) = p_2, \ 
   \psi_2( 1 ) = (-p_2, p_1)
$$
and
$$ p_1 = 2 y h \partial_x+3 x^2 \partial_y ,
   p_2 = 2 x \partial_x+3 y \partial_y+6 h^2.
$$
So, we cannot replace the Schreyer resolution by the minimal
resolution.

Finally, the authors note the iteration approach to get local
cohomology groups as in \cite{Walther} may also improve the performance 
of our algorithms, but it is a future problem.

\bigbreak

During preparing this paper, we knew the paper
Castro-Jim\'enez et al., Homogenizing differential operators, 1997, preprint.
The contents of this section may have an overlap with their paper.

\bigskip
\begin{flushright}
Toshinori Oaku \\
{\tt oaku@math.yokohama-cu.ac.jp} \\
Department of Mathematical Sciences, Yokohama City University \\
Seto 22-2, Kanazawa-ku, Yokohama, 236-0027 Japan
\end{flushright}
\begin{flushright}
Nobuki Takayama\\
{\tt takayama@math.kobe-u.ac.jp}\\
Department of Mathematics, Kobe University \\
Rokko, Kobe, 657-8501 Japan
\end{flushright}
\end{document}